\documentclass[onefignum,onetabnum]{siamart190516}

\usepackage{lipsum}
\usepackage{amsfonts}
\usepackage{graphicx}
\usepackage{epstopdf}
\ifpdf
  \DeclareGraphicsExtensions{.eps,.pdf,.png,.jpg}
\else
  \DeclareGraphicsExtensions{.eps}
\fi

\usepackage{graphicx}
\usepackage{amssymb}

\usepackage{hyperref}
\usepackage{bm}
\usepackage{amsmath}
\usepackage{mathtools}
\usepackage{listings}
\usepackage{color}
\usepackage{siunitx}
\usepackage{booktabs}
\usepackage{algpseudocode}
\usepackage{tikz}
\usetikzlibrary{positioning}
\usetikzlibrary{arrows}
\usetikzlibrary{shapes.multipart}

\newcommand{\timepart}{\mathcal{I}}

\newcommand{\stfesd}{S^{\bf p}(\timepart; \mathcal{D}_h)}

\newcommand{\overbar}[1]{\mkern 1.5mu\overline{\mkern-1.5mu#1\mkern-1.5mu}\mkern 1.5mu}

\newcommand*\diff{\mathop{}\!\mathrm{d}}

\DeclareMathOperator{\Diam}{diam}

\DeclarePairedDelimiter\floor{\lfloor}{\rfloor}

\newcommand{\ud}{\,\mathrm{d}}
\newcommand{\mean}[1]{\{\!\!\{#1\}\!\!\}}                
\newcommand{\jump}[1]{[\![#1]\!]}                        
\newcommand{\ujump}[1]{\lfloor #1\rfloor}   
\newcommand{\fes}{S^{{\bf p}}_{\mathcal{T}_h}}

\newcommand*{\defeq}{\mathrel{\vcenter{\baselineskip0.5ex \lineskiplimit0pt
			\hbox{\scriptsize.}\hbox{\scriptsize.}}}=}

\newcommand\restr[2]{{
		\left.\kern-\nulldelimiterspace 
		#1 
		\vphantom{\big|} 
		\right|_{#2} 
}}

\newcommand{\Rev}[1]{{#1}}

\newsiamremark{remark}{Remark}
\newsiamremark{hypothesis}{Hypothesis}
\crefname{hypothesis}{Hypothesis}{Hypotheses}
\newsiamthm{claim}{Claim}

\headers{GPU-accelerated dG methods on polytopic meshes}{Z.~Dong, E.~H.~Georgoulis, and T.~Kappas}

\title{GPU-accelerated discontinuous Galerkin methods on polytopic meshes\thanks{
\funding{The authors gratefully acknowledge the support by The Leverhulme Trust (grant RPG-2015- 306).  Also, this research work was supported by the Hellenic Foundation for Research and Innovation (H.F.R.I.) under the ``First Call for H.F.R.I. Research Projects to support Faculty members and Researchers and the procurement of high-cost research equipment grant'' (Pro\-ject Number: 3270). }}}

\author{ZHAONAN DONG\thanks{Inria, 2 rue Simone Iff, 75589 Paris, France \&  CERMICS, Ecole des Ponts, 77455 Marne-la-Vall\'{e}e 2, France (\email{zhaonan.dong@inria.fr}) } 
	\and Emmanuil H. Georgoulis\thanks{School of Mathematics and Actuarial Science, University of Leicester, LE1 7RH, United Kingdom (\email{Emmanuil.Georgoulis@le.ac.uk}) \& Department of Mathematics, School of Applied Mathematical and Physical Sciences, National Technical University of Athens, Zografou 15780, Greece \& IACM-FORTH, Heraklion, Crete, Greece}
	\and Thomas Kappas\thanks{School of Mathematics and Actuarial Science, University of Leicester, LE1 7RH, United Kingdom (\email{tk223@le.ac.uk})}
}

\usepackage{amsopn}

\ifpdf
\hypersetup{
  pdftitle={GPU-accelerated discontinuous Galerkin methods on polytopic meshes},
  pdfauthor={Z. Dong, E. H. Georgoulis, and T. Kappas}
}
\fi

\begin{document}

\maketitle


\begin{abstract}
	Discontinuous Galerkin (dG) methods on meshes consisting of polygonal/polyhedral (henceforth, collectively termed as \emph{polytopic}) elements have received considerable attention in recent years. Due to the physical frame basis functions used typically and the quadrature challenges involved, the matrix-assembly step for these methods is often computationally cumbersome. To address this important practical issue, this work proposes two parallel assembly implementation algorithms on CUDA-enabled graphics cards for the interior penalty dG method on polytopic meshes for various classes of linear PDE problems. We are concerned with both single GPU parallelization, as well as with implementation on distributed GPU nodes. The results included showcase almost linear scalability of the quadrature step with respect to the number of GPU-cores used, since no communication is needed for the assembly step. In turn, this can justify the claim that polytopic dG methods can be implemented extremely efficiently, as any assembly computing time overhead compared to finite elements on `standard' simplicial or box-type meshes can be effectively circumvented by the proposed algorithms.
\end{abstract}



\section{Introduction}
Discontinuous Galerkin (dG) methods have received considerable attention during the last two decades. By combining advantages from both finite element methods (FEMs) and finite volume methods (FVMs) they allow the simple treatment of complicated computational geometries, ease of adaptivity and stability for non-self-adjoint PDE problems \cite{hss,MR1842161}. 

More recently, dG approaches have been shown to be applicable on extremely general computational meshes, consisting of general polytopic elements with \Rev{an} arbitrary number of faces and different local elemental polynomial degrees \cite{DGpoly1,MR2846986,DGpoly2,DGpolyparabolic,DGpolybook}. 
A basic feature of these methods is the use of physical frame polynomial bases, as opposed to the standard practice of mapped basis functions in standard finite element implementations. The presence of physical frame basis functions together with the highly involved quadrature requirements over polytopic elements pose new algorithmic complexity challenges in the context of matrix assembly.

The implementation of arbitrary order quadrature rules for non-poly\-nomial integrands over general polytopic domains is highly non-trivial and is addressed in the literature through various techniques. The most general and widely used approach is the subdivision of polytopic elements into basic simplicial or prismatic (with simplicial or hypercubical bases) sub-elements; standard quadrature rules are then employed on each sub-element \cite{DGpoly1,MR2846986,DGpolyparabolic,DGpolybook}. Alternative approaches include the use of Euler's formula for homogeneous polynomials, see, e.g., \cite{gauss_0,gauss_1,gauss_2, gauss_dg}, or the direct derivation of quadrature points for general polytopes, see, e.g., \cite{vianello_2,vianello_3}.

The use of subdivisions, when implemented serially, is typically computationally demanding in this context. At the same time, this approach guarantees the quality of the assembled matrices for nonlinear problems and for problems with localised heterogeneous coefficients. Moreover, in this context quadrature-related variational crimes can be controlled both theoretically and in practice. In contrast, the approaches using Euler's formula \cite{gauss_0,gauss_1,gauss_2, gauss_dg} are typically faster for quadrature computations on polytopic meshes for \emph{polynomial} integrands only; they are not known to offer safeguarded quadrature error control for \emph{non-polynomial} integrands. As such, they cannot be used with confidence in the assembly of nonlinear problems or for problems with non-polynomial PDE coefficients. Finally, tailored quadrature rules for polytopes \cite{vianello_2,vianello_3} typically require costly prior quadrature point optimization steps and, therefore, are not suitable for meshes with highly variable element shapes. Nonetheless, matrix assembly is a \emph{highly parallelizable} process. Thus, it is possible to take advantage of modern computer architectures to achieve highly efficient implementation of all the above approaches.  

Graphics Processing Units (GPUs) have been traditionally used for graphics output to a display device. GPUs offer widely available parallel processing capacity, and they are typically more economical in terms of floating point operations per Watt of electricity \cite{GPU-CPU} than CPU clusters of similar parallelization specifications for basic multiply-and-accumulate processes. At the same time, GPUs can only achieve high performance parallelization on substantially more restricted data structures than CPU clusters. Fortunately, the basic multiply-and-accumulate structure of standard quadrature rules can be naturally implemented within the fast operating ranges of modern GPUs; we refer to \cite{MR3491958} for a GPU-accelerated implementation for low order conforming elements,
or \cite{MR3170331} for a respective high order study, showcasing the acceleration potential in this context
. Also, in \cite{MR3001216} a suite of algorithms and hardware are tested, for low order elements, including usage of atomic operations to avoid race conditions, and mesh `coloring' allowing \Rev{invocation of} different kernels for differently colored patches. 

The benefits of GPU-acceleration in the context of discontinuous Galerkin methods have been studied extensively in the literature over the last decade or so for various classes of electromagnetic, fluid flow and other hyperbolic PDE problems; we refer to \cite{MR2573336,MR3479015,MR3503992,MR3941414,MR3340999,MR3540344,MR3639583} for some of the most successful results in the area. The predominant application setting involves explicit time-stepping, e.g., by structure-preserving Runge-Kutta methods, combined with discontinuous Galerkin spatial discretizations with nodal representation of local finite element spaces \cite{MR2372235}. Indeed, \Rev{owing} in part to the characteristic minimal communication between elemental spaces, \Rev{irrespective} of the local polynomial degree used, dG methods have shown impressive run-time acceleration when implemented in GPU architectures \cite{MR2573336,MR3503992}. 

The present work is concerned with the development and performance study of GPU-accelerated \emph{assembly algorithms} for dG methods on unstructured meshes \Rev{comprising} \emph{extremely general polytopic elements}. This is achieved via a novel CUDA implementation of the $hp$-version interior penalty dG method for equations with non-negative characteristic form on polytopic meshes in $\mathbb{R}^d$. This class of equations, which includes elliptic, parabolic and first order hyperbolic PDEs, as well as equations of changing type, offers a sufficiently general setting for software development. With regard to shape generality, each element is allowed to be a general polytope with \emph{arbitrary} number of $(d-1)$-dimensional polytopic faces; we refer to \cite{DGpolybook} for a detailed discussion on the definition and structure of dG methods on polytopic meshes. The element-shape generality requires both new data structures as well as the resolution of new algorithmic challenges, compared to dG implementations on standard simplicial \Rev{or} box-type meshes \cite{MR2573336,MR3479015,MR3503992,MR3941414,MR3340999,MR3540344,MR3639583}. The algorithms presented below aim use parallelization within GPU clusters to address the key challenge of reducing the computational cost of arbitrary order quadrature rules over general polytopic domains. Given the extreme scalability potential, the quadratures are performed via subdivisions of the polytopic elements into basic simplicial or prismatic sub-elements. \Rev{Standard quadrature rules are, in turn, employed on each of these sub-elements before being accumulated into a matrix entry.}  Correspondingly, for the computation of the face contributions in the present $hp$-version dG setting, subdivision of the $(d-1)$-dimensional faces into simplicial/quadrilateral sub-faces is performed. The choice of quadrature method is made specifically to enable universal applicability: assembly of implicit methods for nonlinear problems or of highly heterogeneous PDE coefficients is possible. We stress, however, that a CUDA implementation of Euler-formula methods, e.g., the one proposed in \cite{gauss_dg}, as by all means possible within the presented algorithmic development.  As we shall see below, the excellent scalability of the implementation essentially removes the computational overhead due to subdivision: quadrature kernels require comparable or less time than the sorting algorithms used to process the resulting arrays. To highlight the performance and the versatility of the proposed algorithms, we consider the interior penalty dG method for:
\begin{enumerate}
	\item[a)]  fully $d$-dimensional ($dD$) problems, $d=2,3$, with non-negative characteristic form, approximated on $dD$ unstructured polytopic meshes, and 
	\item[b)] $(dD+1)$ space-time parabolic problems, $d=2,3$, approximated on prismatic space-time elements with polytopic bases, with the prism bases perpendicular to the time direction. The dG method in this case is equivalent to a combined dG-timestepping scheme with interior penalty dG discretization in space.
\end{enumerate}

For completeness, we present and compare two distinct algorithmic approaches for the matrix assembly: 
\begin{enumerate}
	\item first compute quadrature values for each simplex in the simplicial subdivision; then combine the values appropriately corresponding to each polytopic element;
	
	\item  first precompute the final matrix sparsity patterns, then compute quadrature values and populate the matrices.
\end{enumerate}
The first approach is typical in finite element codes on simplicial or box-type meshes. Interestingly, in the context of general unstructured polytopic meshes with ``many'' faces per element the first approach produces large number of duplicate values in the sparse array formats; these require further costly sorting and processing. On the other hand, the second approach lends itself more naturally to the case of general polytopic elements with arbitrary number of faces per element. The CUDA implementation of the proposed algorithms is able to achieve small run-times for very large discretizations in both $3D$ and $(2D+1)$ settings. Further, we investigate also the scalability of the second approach in a parallel architecture comprising multiple GPUs; the implementation is carried out using basic Message Passing Interface (MPI) tools. \Rev{More specifically, the MPI implementation breaks the problem into completely \emph{independent} processes, each assembling for a part of the mesh using a CUDA-enabled GPU. The actual matrix creation takes place through CUDA-enabled GPUs, with each CUDA thread (the smallest execution unit in a CUDA program,) calculating individual matrix entries.}

The remainder of this work is structured as follows. The general model PDE problem and some important special cases are presented in Section \ref{sec: modelaki}. Section \ref{S:2} contains detailed description of the approaches 1. and 2. above in the implementation of polytopic dG methods on GPUs, along with some initial numerical experiments highlighting the superior performance of the second approach. Finally, in Section \ref{sec:numerics}, we present a number of challenging numerical experiments on single GPU and multiple GPUs.

\section{Model problem and discretization}\label{sec: modelaki}

Let $\Omega$ be a bounded open polygonal/polyhedral domain in $\mathbb{R}^d$, $d=2,3,4$. 

\subsection{Equations with non-negative characteristic form} We consider the advection-diffusion-reaction equation 
\begin{equation}\label{CDR}
	- \nabla \cdot (A \nabla u)+ \bold{b} \cdot \nabla u +cu = f, \quad \mbox{ in } ~ \Omega,
\end{equation}
where $c\in L_\infty (\Omega)$, $f\in L_2 (\Omega)$,  and $\bold{b} :=(b_1,b_2,\dots,b_d)^\top \in [W_{\infty}^1(\Omega)]^d$. Here, $A = \{ a_{ij} \}_{i,j=1}^d$  is a symmetric positive semidefinite tensor whose entries $a_{ij}$ are bounded,
real-valued functions defined on $\bar{\Omega}$, with 
\begin{equation*}\label{diffusion_tensor}
	\bm{\xi}^\top A(x) \bm{\xi} \geq 0 \quad  \forall \bm{\xi} \in  \mathbb{R}^d, \quad \text{a.e.} \quad x\in \bar{\Omega}.
\end{equation*}
Under the above hypothesis, \eqref{CDR} is termed \emph{a partial differential equation with non-negative characteristic form}.

We denote by  $\bold{n}(x) =            \{ n_i(x) \}^d_{i=1}$ the unit outward normal vector to $\partial \Omega$ at $x\in \partial \Omega$ and introduce 
\begin{equation*}\label{Fichera}
	\begin{aligned}
		 \partial \Omega_0 = \big{  \{  } x\in \partial \Omega:  & \bold{n}(x)^\top \Rev{A}(x) \bold{n}(x) >0 \big{  \}  } ,\\
		 \partial_- \Omega = \big{  \{  } x \in \partial \Omega \backslash \partial \Omega_0:  \bold{b}(x) \cdot \bold{n}(x)< 0 \big{  \}  } &,  \quad 
		 \partial_+ \Omega = \big{  \{  } x \in \partial \Omega \backslash \partial \Omega_0: \bold{b}(x) \cdot \bold{n}(x) \geq 0 \big{  \}  } .
	\end{aligned}
\end{equation*}
The sets $\partial_- \Omega$ and $\partial_+ \Omega$ are referred to as the inflow and outflow boundary, respectively. Note that  $\partial \Omega = \partial \Omega_0 \cup \partial_- \Omega  \cup \partial_+ \Omega$.  If $\partial \Omega_0$ is nonempty, we subdivide it into two disjoint subsets $\partial \Omega_{\rm D}$ and $\partial \Omega_{\rm N}$, \Rev{with $\partial\Omega_{\rm D}$ nonempty and relatively open in $\Omega_0$,} on which we consider the boundary conditions:
\begin{equation}\label{boundary_condition}
	u = g_{\rm D} \quad \text{on} \quad \partial \Omega_{\rm D}\cup \partial \Omega_- , \quad \bold{n}\cdot (A \nabla u) = g_{\rm N} \quad \text{on} \quad \partial \Omega_{\rm N},
\end{equation}
and also adopt the hypothesis that $\bold{b} \cdot \bold{n} \geq 0$ on $\partial \Omega_{\rm N}$, whenever $\partial \Omega_{\rm N}$ is nonempty. Additionally, assuming that there exists a positive constant $\gamma_{0}$ such that
$
	c_0(x)^2:=c(x)-1/2\nabla\cdot \bold{b}(x)\geq \gamma_0 $ a.e.  $x\in\Omega$, the well-posedness of the boundary value problem \eqref{CDR}, \eqref{boundary_condition} follows.

\subsection{Discontinuous Galerkin method}

Let $\mathcal{T}_h$ be a subdivision of $\Omega$ into disjoint open \Rev{polygonal or polyhedral elements $\kappa$ ($d = 2$ or $d = 3, 4,$ respectively) such that} \ $\overbar{\Omega}=\cup_{\kappa \in \mathcal{T}_h} \overbar{\kappa}$ and set  $h_{\kappa} \defeq \Diam(\kappa)$. Let also $\mathcal{F}_h$ be the set of all open $(d-1)$--dimensional hyperplanar faces associated with $\mathcal{T}_h$. We write $\mathcal{F}_h = \mathcal{F}_h^{I} \cup \mathcal{F}_h^{B}$, \Rev{with $ \mathcal{F}_h^{B}$ the set of all boundary faces and $\mathcal{F}_h^{I}$ denotes the set of all interior faces, i.e. the faces} shared by two elements. By allowing general polytopic elements in $\mathcal{T}_h$, it is by all means possible that two elements share more than one \Rev{face}. Nonetheless, the term \emph{face} will refer to a $(d-1)$-dimensional planar region of each element, while the term \emph{interface} will refer to the totality of the common boundary between two elements. The domain of all (interior) interfaces will be denoted by $\Gamma_{\rm int}:=\cup_{F\in\mathcal{F}_h^I}F\subset\Omega$.


Given $\kappa \in  \mathcal{T}_h$, we write $p_{\kappa}\in\mathbb{N}$ to denote the \emph{polynomial degree} of the element $\kappa$, and collect the $p_{\kappa}$ in
the vector ${\bf p}:=(p_{\kappa}:\kappa \in  \mathcal{T}_h)$. We then define the \emph{finite element space} $\fes$ with
respect to $\mathcal{T}_h$ and ${\bf p}$ by
\[
\fes:=\{u\in L_2(\Omega)
:u|_{\kappa}\in\mathcal{P}_{p_{\kappa}}(\kappa),\kappa \in  \mathcal{T}_h\},
\]
where
$\mathcal{P}_{p_\kappa}(\kappa)$ denotes the space of polynomials of total degree $p_\kappa$ on $\kappa$.  Note  that the local elemental polynomial spaces employed within the definition of $\fes$ are defined in the physical coordinate system, without the need to map from a given reference or canonical frame.

Next, we introduce some trace \Rev{operators} used in the definition of discontinuous Galerkin methods.  For element $\kappa\in \mathcal{T}_h$, we define the \emph{inflow} and \emph{outflow} parts of its boundary $\partial \kappa$ by
\begin{eqnarray*}
	\partial_{-}\kappa=\{x\in \partial\kappa, \quad \bold{b}(x)\cdot{\bold{n}_\kappa(x)}<0\} ,  \quad \partial_{+}\kappa=\{x\in \partial\kappa, \quad \bold{b}(x)\cdot{\bold{n}_\kappa(x)} \geq 0\},
\end{eqnarray*}
respectively, with $\bold{n}_\kappa(x)$ denoting the unit outward normal vector to $\partial \kappa$ 
at $x \in \partial\kappa$. 

We shall also make use of the \emph{upwind jump} of a function $v$ across a face $F \subset \partial_-\kappa \backslash \partial \Omega$, denoted by
\begin{equation*}
	\ujump{v}:=v^+_\kappa-v^-_\kappa .
\end{equation*}
Also, for $\kappa_i, \kappa_j\in {\cal T}_h$ two 
adjacent elements sharing a face $F =\partial \kappa_i \cap \partial \kappa_j\subset \mathcal{F}_h^{I}$, we write $\bold{n}_{i}$ and $\bold{n}_{j}$  to denote the 
outward unit normal vectors on $F$, relative to $\kappa_i$ and $\kappa_j$, respectively. 
\Rev{Let $w$  be  a (scalar- or vector-valued) function and $v$ be a scalar  function, } that are smooth enough on each element to have a well-defined trace on $F$ from within both $\partial \kappa_i, \partial \kappa_j$. \Rev{We define the \emph{average} of $w$ and \emph{jump} of $v$ across $F$ by
$$
\mean{w}_F:=\frac{1}{2}(w|_{F\cap\partial \kappa_i}+w|_{F\cap\partial \kappa_j}),\quad \jump{v}_F :=v|_{F\cap\partial\kappa_i}\cdot{\bold{n}}_{i}+v|_{F\cap\partial \kappa_j}\cdot{\bold{n}}_{j}, $$
respectively. On a boundary face $F\subset  \mathcal{F}_h^{B}$, with $F \subset \partial \kappa_i$,
we simply set 
$\mean{w}_F=w|_{F\cap\partial \kappa_i}$,  $  \jump{v}_F :=v|_{F\cap\partial \kappa_i}\cdot{\bold{n}}_{i},
$}
noting that in the last case~$\bold{n}_{i}$ coincides with the unit outward normal
vector on the boundary $\partial \Omega$ if the domain is represented exactly by the mesh.  We observe that $\jump{\cdot }$ and $\ujump{\cdot}$ may differ only up to a sign. For brevity, we also define the \emph{broken gradient} $\nabla_h v$ of a sufficiently smooth function $v$ to be given by $\nabla_h v|_\kappa=(\nabla v)|_\kappa$ for all $\kappa\in\mathcal{T}_h$.

The \Rev{symmetric} interior penalty discontinuous Galerkin method for  \eqref{CDR}, \eqref{boundary_condition} is given by: find $u_h\in \fes$ such that
\begin{equation}\label{adv_galerkin_dg}
	B(u_h,v_h)=\ell(v_h) ,\quad\text{ for all \ }v_h\in \fes;
\end{equation}
the bilinear form $B(\cdot, \cdot ):\fes\times \fes\to \mathbb{R}$ is given by:
\begin{equation}\label{dg-bilinear}
	\begin{aligned}
		&B(w,v) :=  \int_{\Omega} \big( A \nabla_h w \cdot \nabla_h v+  \bold{b} \cdot  \nabla_h w v  +cw v \big) \ud x  \\
		&-\sum_{\kappa \in \mathcal{T}_h} 
		\int_{\partial_-\kappa \backslash \partial \Omega}(\bold{b} \cdot  \bold{n})\ujump{w}v^+\ud s
		-\sum_{\kappa \in \mathcal{T}_h} \int_{\partial_-\kappa \cap (\partial \Omega_{\rm D} \cup \partial \Omega_-)}(\bold{b} \cdot \bold{n})w^+v^+\ud s \\
		& -
		\int_{\Gamma_{\rm int}\cup \partial\Omega_{\rm D}} \big( \mean{A \nabla w} \cdot \jump{v} + \mean{A \nabla v} \cdot \jump{w} - \sigma \jump{w} \cdot \jump{v} \big) \ud s,
	\end{aligned}
\end{equation}
and the linear functional $\ell:\fes\to\mathbb{R}$ by
\begin{equation}\label{dg-linear}
	\begin{aligned}
		\ell(v)
		:= &\int_{\Omega} f  v \ud x
		-\sum_{\kappa \in \mathcal{T}_h} \int_{\partial_-\kappa \cap (\partial \Omega_{\rm D} \cup \partial \Omega_-)}(\bold{b} \cdot \bold{n}) g_{\rm D}v^+\ud s  \\
		&- \int_{\partial\Omega_{\rm D}} g_{\rm D}^{} \Big( (A \nabla v ) \cdot \bold{n} -\sigma v \Big)\ud s  
		+ \int_{\partial\Omega_{\rm N}} g_{\rm N}^{}v \ud s. 
	\end{aligned}
\end{equation}
The nonnegative function  $\sigma : \Gamma_{\rm int}\cup \partial \Omega_{\rm D}\rightarrow \mathbb{R}$ appearing in~\eqref{dg-bilinear} and~\eqref{dg-linear} is referred to as the \emph{discontinuity-penalization parameter}; its precise definition used in the numerical experiments below is provided in Appendix \ref{penalty}.

\subsection{An important special case}

In this section, we focus on the important subclass of first order (in time) evolution problems which can be considered as a special case of PDEs with non-negative characteristic form. In particular, for $d=3,4$, we consider the following special form of $A:\mathbb{R}^d \to \mathbb{R}^{d\times d}$ and ${\bf b}: \mathbb{R}^d \rightarrow \mathbb{R}^{d}$, 
\[
A = \bigg(\begin{array}{cc}
{\bf a}\ & {\bf 0} \\ 
{\bf 0}^\top\ & 0
\end{array} \bigg),\quad\qquad {\bf b}= \bigg( \begin{array}{c}
{\bf {\bf w}}\\
1
\end{array} \bigg),
\]
with ${\bf a}:\mathbb{R}^{s} \to \mathbb{R}^{s\times s }$, $s:=d-1$, a symmetric non-negative definite tensor and ${\bf w}:\mathbb{R}^{s} \to \mathbb{R}^{s }$ the spatial wind/advection direction.  Substituting this selection into \eqref{CDR} gives rise to the classical first order (in time) evolution equations, with the last variable designating the time direction, viz.,
\begin{equation}\label{first_order_evol}
	\partial_t u  - \nabla \cdot ({\bf a} \nabla u) +{\bf w}\cdot \nabla u+cu= f  \quad \text{in } \Omega , 
\end{equation}
with $\nabla \defeq (\partial_{x_1},  \dots,\partial_{ x_s})^\top$ and $\nabla \cdot \defeq \sum_{i=1}^s\partial_{x_s}$  the gradient and divergence operators with respect to the spatial variables only, respectively. If ${\bf{a}} $ is additionally uniformly positive definite, i.e.,  
\begin{equation}\label{spd}
	\bm{\zeta}^{\top} {\bf{a}}(t,x) \bm{\zeta} 
	\geq \theta |\bm{\zeta}|^{2}
	> 0 \quad \forall \, \bm{\zeta} \in \mathbb{R}^{s}, \quad \textnormal{a.e.  in } \Omega ,
\end{equation}
with \Rev{$\theta$ a positive constant and}  $x:=(x_1,\dots,x_s)^\top$, \eqref{first_order_evol} is, in particular, a parabolic PDE. The preference on temporally implicit high order discretizations for parabolic PDE problems motivates the use space-time dG methods for these problems. Moreover, if spatio-temporal variability is present in the PDE coefficients, matrix assembly has to be performed for every/most time-steps. Hence, it is of interest in the present context to focus also in space-time methods for parabolic problems. To that end, we assume that \eqref{spd} holds for the rest of this section.

Let $\Omega  \defeq J\times D$ be a space-time domain with $\Rev{J\defeq (0,T]} \subset \mathbb{R}^+$ a time interval, and $D\subset \mathbb{R}^{s}$ denoting the spatial domain.  Let also $\partial D_{\rm D}$ denote the Dirichlet part of the boundary $D$.  Finally, we set $\partial D_{\rm N} :=\partial D\backslash\partial D_{\rm D}$ for the Neumann boundary. Thus, we have $\partial\Omega_{\rm D}= \{(t,x): t\in J, x\in \partial D_{\rm D}(x)\}$, and correspondingly for $\partial\Omega_{\rm N}$. 

We consider the linear parabolic problem:
\begin{align}\label{parabolic-problem1}
	\partial_t u  - \nabla \cdot ({\bf a} \nabla u) +{\bf w}\cdot \nabla u+cu=&\ f  \ \quad \text{in }\ \Omega=J\times D , 
	\nonumber\\
	u=&\ u_0 \quad \text{on }\ \{0\}\times D , \\
	u =&\ g_{\rm D}  \quad \text{on }\ J \times \partial D_{\rm D}, \nonumber\\
	 \bold{n}\cdot (\bold{a} \nabla u)  =&\ g_{\rm N}  \quad  \text{on }\ J \times \partial D_{\rm N}, \nonumber 
\end{align}
with  $f \in L_2(J;L_2(D))$, ${\bf{a}} \in [L_{\infty}(\Omega)]^{s\times s }$, ${\bf w}\in [W^1_\infty(\Omega)]^s$, ${c} \in L_{\infty}(\Omega)$, $u_0 \in L_2(D)$; \Rev{the initial condition at $t=0$},  $g_{\rm D}\in L_2 (J; H^{1/2}( \partial D_{\rm D} ))$ and $g_{\rm N}\in L_2 (J; H^{-1/2}( \partial D_{\rm N} ))$; the Dirichlet and/or Neumann boundary conditions, which may be time-varying, are imposed only on $J \times \partial D$, since the boundary portions $\{0,T\}\times D\subset \partial\Omega\backslash \partial\Omega_0$.

Although it is by all means possible to consider an unstructured space-time mesh for this problem also, we prefer to use the structure of the equation and construct a mesh based on space-time slabs. That way, it is possible to solve for each slab independently of the next ones in the time direction; this idea goes back to \cite{jamet}. 

To that end, we begin by introducing a temporal discretization first. Let $\timepart:=\{ I_n \}_{n=1}^{N_t}$ be a partition of the time interval $J$ into $N_t$ time steps $I_n$, with $I_n=(t_{n-1}, t_n]$, \Rev{for a set of time nodes $\{ t_n \}_{n=0}^{N_t}$ with $0 = t_0 < t_1 < \ldots < t_{N_t} = T$. Let also  $\tau_n := t_n-t_{n-1}$ denote} the length of $I_n$.  Let also $\mathcal{D}_h$ an $s$-dimensional polytopic mesh subdividing the spatial domain $D$. Each space-time slab $I_n\times D$, is then subdivided into a mesh $\mathcal{T}_h$ comprising disjoint open prismatic elements $\kappa_n:=I_n\times\kappa$, $\kappa\in \mathcal{D}_h$; thus, we have $\mathcal{T}_h:=\timepart\times \mathcal{D}_h$.

For notational simplicity we do not include explicitly `local' timestepping within one slab in the discussion, i.e., elements $\kappa_n^i:=I_n^i\times\kappa$ arising by subdividing $I_n$ into disjoint subintervals $I_n^i$ for some $i=1,\dots, m_\kappa$, \Rev{with $m_\kappa$ being  the maximum number of local stepping on $\kappa_n^i$}. We stress, however, that this capability is included in the computer implementation presented below. Also, in the present space-time slabbing setting, we assume that the Dirichlet or Neumann \emph{domains} are time-independent; if they are required to be so, we can revert back to the `monolithic' formulation \eqref{adv_galerkin_dg}. Still, the Dirichlet or Neumann data $g_{\rm D}^{}$ and $g_{\rm N}^{}$, respectively, are allowed to be time-dependent functions.

Let $p_{\kappa_n}$ denote the polynomial degree of the space--time element $\kappa_n$ and let $\bm{p}_n$ be a vector of the polynomial degrees of all elements in $I_n\times \mathcal{D}_h$. We define the space-time finite element space for the time interval $I_n$ by
\begin{equation*}
	V^{\bm{p}_n}(I_n;\mathcal{D}_h) \defeq \{ u \in L^2(I_n \times D): \restr{u}{\kappa_n} \in \mathcal{R}_{p_{\kappa_n}}(\kappa_n), \kappa_n \in I_n \times  \mathcal{D}_h \},
\end{equation*}
where $\mathcal{R}\in\{\mathcal{P}, \mathcal{PQ} \}$; $\mathcal{P}_{p_{\kappa_n}}(\kappa_n)$ denotes the space of polynomials of space-time total degree $p_{\kappa_n}$ on $\kappa_n$; $\mathcal{PQ}_{p_{\kappa_n}}(\kappa_n)$ \Rev{is} the space of polynomials of degree $p_{\kappa_n}$ in temporal variable tensor-product with polynomials of total degree $p_{\kappa_n}$ in all spatial variables. $\mathcal{PQ}_{p_{\kappa_n}}$ is the standard choice when $\kappa$ are simplices \cite{johnson_thomee}, while $\mathcal{P}_{p_{\kappa_n}}$ was proposed in \cite{DGpolyparabolic} as a reduced complexity choice, without loss of rate of convergence in energy-like norms.  The global space-time finite element space $\stfesd$ is then defined by
\begin{equation*}
	\stfesd\defeq \bigoplus_{n=1}^{N_t} V^{\bm{p}_n}(I_n;\mathcal{D}_h),
\end{equation*}
with $\bm{p}:=(\bm{p}_1,\bm{p}_2,\dots,\bm{p}_n)$ the array of elemental polynomial degrees.
Note that the local elemental space-time polynomial spaces are defined in the physical coordinate system, i.e., without mapping from a reference frame; this ensures optimal approximation capabilities over each space-time prism. 
Finally, let $h_{\kappa_n}$ denote the diameter of each space-time element $\kappa_n$ and let $\bm{h} \defeq (h_{\kappa_n}: \kappa_n \in \stfesd)$ the array of elemental diameters.

The specific space-time slab mesh structure (which is orthogonal with respect to the time direction) allows for reorganization of the face terms in \eqref{dg-bilinear} and \eqref{dg-linear}. In particular, we set $\Gamma_{\rm int}^{n}:= I_n \times \Gamma_{\rm int}$, i.e., the union of all element faces in the space-time slab $I_n\times D$ parallel to the time direction. \Rev{Additional} to the operators $\mean{\cdot}$ $\jump{\cdot}$, $\ujump{\cdot}$, we also specify the \emph{time-jump}  $\ujump{u}_n:=u_n^+-u_n^-$, where $u_n^{\pm} \defeq \lim_{s \to 0^{\pm}} u(t_n+s)$, $1 \leq n \leq N_t$ and $u_0^+:=\lim_{s \to 0^+} u(s)$.
We note that the time-jump is just \Rev{an} instance of the upwind jump $\ujump{\cdot}$ operator across elemental faces perpendicular to the time direction, i.e., the prismatic element bases. For brevity, we also denote by $(\cdot,\cdot)_{\omega}$ the $L^2(\omega)$-inner product for any measurable set $\omega\subset \mathbb{R}^{d}$. 

The (space--time) discontinuous Galerkin method for \eqref{parabolic-problem1} then becomes: find $u_h \in \stfesd$ such that:
\begin{equation}\label{spacetime_dg}
	B(u_h,v_h) = \ell (v_h) \quad \text{ for all } \, v_h \in \stfesd,
\end{equation}
where
\begin{equation*}
	B(u,v) = \sum_{n=1}^{N_t} \int_{I_n} \big( (\partial_t u, v)_D + B_{\mathrm{S}}(u,v) \big) \diff t + \sum_{n=2}^{N_t}  \left( \floor*{u}_{n-1}, v_{n-1}^{+} \right)_D + (u_0^{+}, v_0^{+})_D, 
\end{equation*}
with the spatial dG bilinear form given by
\begin{equation*}
	\begin{aligned}
		B_{\mathrm{S}}(u,v) :=&  \int_{D} \big( \bold{a} \nabla_h u \cdot \nabla_h v+  \bold{w} \cdot  \nabla_h u v  +cu v \big) \ud x  \\
		&-\sum_{\kappa \in \mathcal{D}_h} \bigg(
		\int_{\partial_-\kappa \backslash \partial D}(\bold{w} \cdot  \bold{n})\ujump{u}v^+\ud s
		+ \int_{\partial_-\kappa \cap \partial D_{\rm D} }(\bold{w} \cdot \bold{n})u^+v^+\ud s\bigg) \\
		& -
		\int_{\Gamma_{\rm int}^n\cup \partial D_{\rm D}} \big( \mean{\bold{a} \nabla u} \cdot \jump{v} + \mean{\bold{a} \nabla v} \cdot \jump{u} - \sigma \jump{u} \cdot \jump{v} \big) \ud s,
	\end{aligned}
\end{equation*}
and the linear functional $\ell $ now becomes
\begin{equation*}
	\ell(v) \defeq \sum_{n=1}^{N_t} \int_{I_n} \bigg( (f,v)_D - \int_{\partial D_{\rm D}} g_{\mathrm{D}} \big( ({\bf{a}} \nabla v)\cdot \bm{n} - \sigma v \big) \diff s 
	+ \int_{\partial D_{\rm N}} g_{\mathrm{N}} v  \diff s
	\bigg) \diff t + (u_0, v_0^+)_D,
\end{equation*}
where $\sigma :J\times (\partial D_{\rm D}\cup \Gamma_{\rm int}) \rightarrow \mathbb{R}_+$ the discontinuity--penalization parameter; see Appendix \ref{penalty} for its definition.

The space-time slabs for this problem allows \Rev{us} to solve for each time-step separately. Thus, \eqref{spacetime_dg} can be solved on each time interval $I_n \in \timepart$, $n=1, \dots ,N_t$.  \Rev{The} solution $U=\restr{u_h}{I_n} \in V^{\bm{p}}(I_n;\mathcal{D}_h)$ is given by:
\begin{multline*}
	\Rev{\int_{I_n} \big( (\partial_t U, W)_D + B_{\mathrm{S}}(U,W) \big) \diff t + (U_{n-1}^{+}, W_{n-1}^{+})_D } \\
	\!= \!\! \Rev{ \int_{I_n} 
	\hspace{-0.2cm}\left( (f,W)_D 
	- \hspace{-0.1cm} \int_{\partial D_{\rm D}} 
	 \hspace{-0.2cm} g_{\mathrm{D}} \big( ({\bf{a}} \nabla W)\cdot \bm {n} - \sigma W \big) \diff s 
	+ \hspace{-0.1cm} \int_{\partial D_{\rm N}} 
	\hspace{-0.2cm}  g_{\mathrm{N}} W  \diff s
	\right) \diff t + (U_{n-1}^{-}, W_{n-1}^+)_D,}
\end{multline*}
for all \Rev{$W \in V^{\bm{p}}(I_n;\mathcal{D}_h)$}, with $U_{n-1}^{-}:=U(t_{n-1}^-)$ the computed solution from the previous time step, $n=2, \dots, N_T$. \Rev{For} $n=1$,  $U_{0}^{-}$ is given by  the initial condition $u_0$.

\section{Implementation on Graphics Processing Units}
\label{S:2}

Quadrature computations during the assembly of the stiffness and mass matrices are computationally demanding steps, especially for high order local polynomial spaces. \Rev{Nonetheless, fast methods for stiffness matrix assembly for standard simplicial or box-type element meshes are widely known, especially for low order elements, resulting to (near) optimal scaling. This is \emph{not} the case, however, to the best of our knowledge, for general polytopic element shapes, especially ones arising by aggressive agglomeration of finer simplicial meshes and the use of physically-defined basis functions. Such a development is deemed important to unlock the complexity reduction potential of polytopic mesh methods, by reducing the cost of assembly, thereby re-instating the linear/nonlinear solvers as the dominant source of complexity on par with classical Galerkin approaches. Therefore, fast quadrature evaluations are highly desirable for general polytopic meshes.}

\Rev{The computational cost of assembly becomes extremely relevant when discretizing evolution PDE problems. Then,} assembly takes place after each time-step when the PDE coefficients and/or boundary conditions depend on the time variable also. In the latter case, in particular, assembly is typically at least as time-consuming as solving for the approximate solution of each time-step. This is because the resulting algebraic systems arising are either block-diagonal (e.g., for hyperbolic problems or when using explicit schemes for parabolic problems), or the convergence of iterative solvers for implicit computations is typically fast: computer solutions from the previous time-step are excellent starting values typically for the iteration. \Rev{It is, therefore, highly relevant to develop algorithms able to perform fast assembly for highly agglomerated polytopic element methods.}

We will describe two approaches taken to accelerate the quadrature computations for dG methods on general polytopic meshes by the use of Graphics Processing Units (GPUs). As discussed above, our aim is to provide an as general as possible framework, allowing in particular to assemble problems with variable PDE coefficients, stemming from, e.g., non-linear or highly heterogeneous equations. 

To that end, the algorithm first subdivides each polytopic element into simplicial sub-elements or, respectively, prismatic ones with simplicial bases in the context of space-time slabs. Correspondingly, for the computation of the face contributions, the algorithm subdivides any co-hyperplanar $(d-1)$-dimensional faces into simplicial sub-faces, or, respectively, rectangular faces for prismatic elements.  For brevity, we shall, henceforth, refer collectively as \emph{simplicial subdivision} of an element for both cases of simplicial and prismatic elements with simplicial faces. Standard quadrature rules are then applied on each sub-element or sub-face. 

We stress, however, that for the special case of polynomial integrands on polytopic domains, it is possible to further accelerate the quadrature computations within the algorithmic framework presented below, via the use of Euler-formula type quadrature rules \cite{gauss_dg}. Nonetheless, we prefer to take the point of view of accelerating the computation for polytopic elements in ``worst-case" scenarios of practical interest, e.g., heterogeneous coefficients and/or nonlinear problems. Moreover, we are also concerned with developing quadratures for extreme case scenarios of polytopic elements with \emph{many} non-co-hyperplanar polygonal faces each, e.g., elements arising via agglomeration of a finer simplicial subdivision. In these cases, we expect that the majority of face integrals will be performed predominantly over triangular faces, thereby the use of Euler-type quadratures is not expected to be significantly advantageous even in the cases they are applicable.  

We note that polytopic elements arising from agglomeration of \emph{many} \Rev{(e.g., tens or, even, hundreds of)} simplicial elements 
of an underlying finer simplicial mesh can be subdivided in a different fashion to the original constituent simplicial elements for the purposes of quadrature. This, in turn allows to minimize the number of sub-elements on which quadratures will be performed. We now discuss the particulars of the parallel calculation of the element and face integrals on a GPU.

\subsection{Computing the integrals}
\begin{figure}[h]
	\begin{minipage}{7cm}
		\centering
		\begin{tikzpicture}[every node/.style={draw,shape=circle,fill=black,inner sep=0pt,minimum size=3pt}]
		
		\node[draw=none,fill=white,above right] at (-1.9,-1)  {\large $\kappa_1$};
		\node[draw=none,fill=white,above left] at (-1,-3)  {\large $\kappa_2$};
		
		\path
		(-3.6,-2.0) node (p0) {}
		(-2.0,-1.8) node (p1) {}
		(0.0,-2.0) node (p2) {}
		(1.5,-1.5) node (p3) {}
		(-2.6,-0.1) node (p4) {}
		(0.1,-0.2) node (p5) {}
		(-3.1,-3.4) node (p6) {}
		(-1.0,-3.6) node (p7) {}
		(1.2,-3.3) node (p8) {};
		
		\draw[line width=0.60mm,black]
		(p0) -- (p1)
		(p1) -- (p2)
		(p2) -- (p3);
		
		\draw[line width=0.30mm,yellow]
		(p0) -- (p1)
		(p1) -- (p2)
		(p2) -- (p3);
		
		\draw[line width=0.30mm,black]
		(p0) -- (p4)
		(p4) -- (p5)
		(p5) -- (p3)
		(p0) -- (p6)
		(p6) -- (p7)
		(p7) -- (p8)
		(p8) -- (p3);
		
		\draw[line width=0.15mm,dashed,black]
		(p1) -- (p4)
		(p1) -- (p5)
		(p2) -- (p5)
		(p1) -- (p6)
		(p1) -- (p7)
		(p2) -- (p7)
		(p2) -- (p8);
		\end{tikzpicture}
	\end{minipage}
	\begin{minipage}{46mm}
		\centering
		\begin{tikzpicture}
		\draw [line width=0.30mm,black] (-1.7,1.7) rectangle (1.7,-1.7);
		\draw [line width=0.15mm,dashed,black,fill=yellow] (-1.4,1.4) rectangle (-1.0, 1.0);
		\draw [line width=0.15mm,dashed,black,fill=yellow] (0.2,-0.2) rectangle (1.0, -1.0);
		\draw [line width=0.15mm,dashed,black,fill=yellow] (0.2,1.4) rectangle (1.0, 1.0);
		\draw [line width=0.15mm,dashed,black,fill=yellow] (-1.4,-0.2) rectangle (-1.0, -1.0);
		\end{tikzpicture}
	\end{minipage}
	\caption{\emph{Left:} Polygonal elements $\kappa_1$, $\kappa_2$ sharing $3$ interior faces. \emph{Right:} Blocks receiving contribution from the shared interior faces.} 
	\label{fig:blocksubmatrices}
\end{figure}
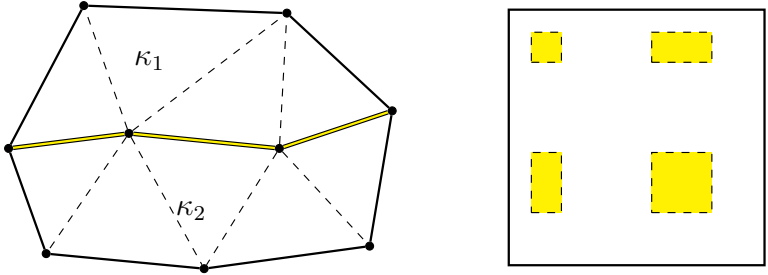

We begin by describing the computation kernels. We implement five kernels: one for volume integrals over each element, and one for each of the face contributions: $\partial\Omega_{\rm D}$ (or $\partial D_{\rm D}$), $\partial\Omega_-$, $\partial\Omega_+$ and $\Gamma_{\rm int}$ (or $\Gamma_{\rm int}^n$), respectively. Implementing different kernels for integrals along boundary and interior faces is recommended, as each process requires a different computation load: for instance the kernel computing interior face integrals on $\Gamma_{\rm int}$ is more computationally demanding since it contributes to four blocks for each interface; we refer to Figure \ref{fig:blocksubmatrices} for an illustration.

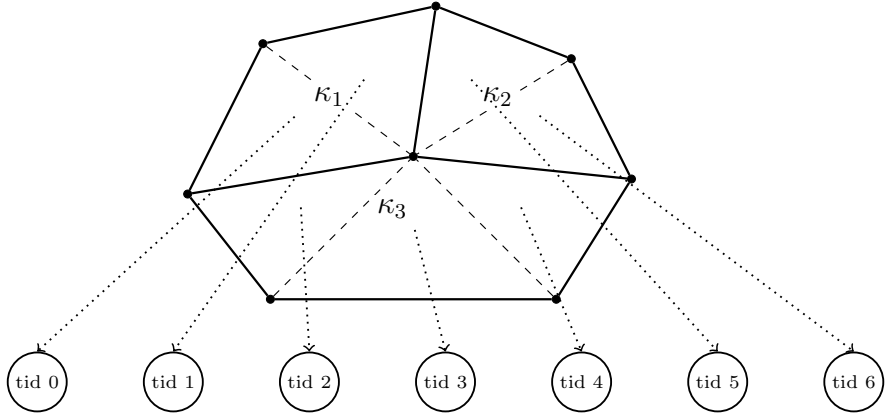
\begin{figure}[!t]
	\centering
	\begin{tikzpicture}[sharp corners=2pt,inner sep=2pt,node distance=.5cm,every text node part/.style={align=center}]
	
	\path[every node/.style={draw,shape=circle,fill=black,inner sep=0pt,minimum size=3pt}]
	(-3,-0.5) node (p0) {}
	(-2,1.5) node (p1) {}
	(0.3,2) node (p2) {}
	(2.1,1.3) node (p3) {}
	(2.9,-0.3) node (p4) {}
	(1.9,-1.9) node (p5) {}
	(-1.9,-1.9) node (p6) {}
	(0,0) node (p7) {};
	
	\draw[line width=0.30mm,black]
	(p0) -- (p1)
	(p1) -- (p2)
	(p2) -- (p3)
	(p3) -- (p4)
	(p4) -- (p5)
	(p5) -- (p6)
	(p6) -- (p0)
	(p0) -- (p7)
	(p2) -- (p7)
	(p4) -- (p7);
	
	\draw[line width=0.15mm,dashed,black]
	(p1) -- (p7)
	(p3) -- (p7)
	(p5) -- (p7)
	(p6) -- (p7);
	
	\node[line width=0.25mm, circle, draw, minimum height = 2mm, minimum width = 2mm] at (-5,-3) (tid0){\scriptsize tid 0};
	\node[line width=0.25mm, circle, draw, right=1cm of tid0, minimum height = 2mm, minimum width = 2mm] (tid1){\scriptsize tid 1};
	\node[line width=0.25mm, circle, draw, right=1cm of tid1, minimum height = 2mm, minimum width = 2mm] (tid2){\scriptsize tid 2};
	\node[line width=0.25mm, circle, draw, right=1cm of tid2, minimum height = 2mm, minimum width = 2mm] (tid3){\scriptsize tid 3};
	\node[line width=0.25mm, circle, draw, right=1cm of tid3, minimum height = 2mm, minimum width = 2mm] (tid4){\scriptsize tid 4};
	\node[line width=0.25mm, circle, draw, right=1cm of tid4, minimum height = 2mm, minimum width = 2mm] (tid5){\scriptsize tid 5};
	\node[line width=0.25mm, circle, draw, right=1cm of tid5, minimum height = 2mm, minimum width = 2mm] (tid6){\scriptsize tid 6};
	
	\node[draw=none,fill=white,below right] at (-1.4,1.0) {\large $\kappa_1$};
	\node[draw=none,fill=white,below left] at (1.4,1.0) {\large $\kappa_2$};
	\node[draw=none,fill=white,below left] at (0,-0.5) {\large $\kappa_3$};
	
	\path
	(-1.5,0.6) node (e0) {}
	(-0.6,1.1) node (e1) {}
	(-1.5,-0.6) node (e2) {}
	(0,-0.9) node (e3) {}
	(1.4,-0.6) node (e4) {}
	(0.7,1.1) node (e5) {}
	(1.6,0.6) node (e6) {};
	
	\draw[->,dotted,line width=0.25mm] (e0) -- (tid0.north);
	\draw[->,dotted,line width=0.25mm] (e1) -- (tid1.north);
	\draw[->,dotted,line width=0.25mm] (e2) -- (tid2.north);
	\draw[->,dotted,line width=0.25mm] (e3) -- (tid3.north);
	\draw[->,dotted,line width=0.25mm] (e4) -- (tid4.north);
	\draw[->,dotted,line width=0.25mm] (e5) -- (tid5.north);
	\draw[->,dotted,line width=0.25mm] (e6) -- (tid6.north);
	
	\end{tikzpicture}
	\caption{Three polygonal elements split into $7$ simplices. One thread is used per simplex and correspondingly for the face integrals. \texttt{tid} refers to ``thread index.''}
	\label{fig:workloadsplit}
\end{figure}

We now discuss the structure of the element and interior face integral computations; the boundary integrals are computed in an analogous, yet simpler, fashion. We allocate one thread for the calculation of the volume integral for each simplex of the simplicial subdivision, as illustrated in Figure \ref{fig:workloadsplit}.
The implementation allows for locally variable polynomial degree. As such, the algorithm groups simplicial subdivision elements with the same polynomial degree, spawning respective number of threads and invoking the element kernel for each polynomial degree. 
This allows for the minimisation of idle threads within each block. For the face integral computations, we are presented with the challenge that two elements sharing the same face may admit different degree basis. To address this, we first group the faces according to the maximum of those two degrees and invoke the interior face integral kernels for each maximum degree separately. The use of hierarchical basis functions allows a unified implementation of the kernels for each degree.

In the standard case of simplicial or mapped box-type elements, the use of nodal basis functions is common practice as it offers significant computational savings \cite{MR2372235,MR2573336,MR3479015,MR3503992,MR3941414,MR3340999,MR3540344,MR3639583}. This is because nodal basis functions allow for offline precomputation of the evaluation of the basis functions on the quadrature points, which are transferred into the physical domain via elemental maps. For general polytopic elements, nodal basis functions are not an option, at least for elements with many vertices. Instead, we employ a different approach: physical domain polynomial basis functions are defined on a rectangular bounding box of each element before being restricted to the polytope \cite{DGpoly1,DGpolybook}. The basis functions of choice in our implementation are tensorized orthonormal Legendre polynomials of \emph{total} degree $p_\kappa$, $\kappa\in\mathcal{T}_h$ (or $p_{\kappa_n}$, respectively). As a result, the innermost \texttt{for-loop} \Rev{in} Algorithm \ref{alg:AssemblyAlgorithm} becomes more expensive, as the evaluation (function \emph{evaluate} in Algorithm  \ref{alg:AssemblyAlgorithm}) of the basis functions on each quadrature point must be calculated in the physical domain directly. 

\begin{algorithm}[h!]
	\begin{algorithmic}[1]
		\State read geometry data
		\For{$i=1$ to $\# Trial \ basis \ functions$}
		\For{$j=1$ to $\# Test \ basis \ functions$}
		\State $val_{ij} \leftarrow 0$
		\For{$k=1$ to $\# Quadrature \ points$}
		\State $u_i = evaluate(i,q_k)$
		\State $v_j = evaluate(j,q_k)$
		\State $val_{ij} \leftarrow val_{ij} + \mathcal{B}(u_i, v_j)$
		\EndFor
		\State write $val_{ij}$ into global memory
		\EndFor
		\EndFor
	\end{algorithmic}
	\caption{Pseudocode for the calculation of the bilinear form for one simplex. 
	}
	\label{alg:AssemblyAlgorithm}
\end{algorithm}

The \emph{evaluation} function first affinely maps the quadrature points from a reference simplex to each sub-simplex of the simplicial subdivision. Then, using a fast and stable polynomial evaluation method, 
each polynomial basis function is evaluated at the respective quadrature points. The polynomial evaluation implementation benefits from the fused multiply--add (FMA) operation capability of the GPU \cite{P754:2008:ISF,DBLP:conf/arith/BoldoM05}. Moreover, storing polynomial coefficients and reference simplex quadrature points in the GPU's \emph{constant memory} can improve memory throughput, as threads belonging to the same warp (a group of 32 threads with contiguous indices) are all executing the polynomial evaluation for the same data set in each loop. Constant memory size of current GPU chips (typically 64KB at the time of writing) is sufficient to hold all the data needed for implementation of local polynomial degree up to $p=10$ in three dimensions.

\subsection{Connectivity}

Upon deciding on a simplicial subdivision for the polytopic mesh, we discuss different potential choices for index and quadrature value allocation manipulations. 
We consider two different \emph{approaches} for the assembly of the stiffness and mass matrices in terms of index manipulation:
\begin{enumerate}
	\item first compute the volume and face quadrature values for each simplex in the simplicial subdivision; then combine the values appropriately corresponding to each polytopic element, or
	
	\item  first precompute the sparsity pattern of the final stiffness and mass matrices, then compute the volume and face quadrature values for each simplex in the simplicial subdivision and, subsequently, populate the matrices with the calculated values.
\end{enumerate}
As we shall see below, each of the above general approaches offer advantages on different settings.

\subsection{Approach 1}

The general idea of the first approach is to calculate the quadrature values for each simplex of the simplicial subdivision, along with their row and column indices in an unsorted coordinate format containing, typically, duplicates. Subsequently, we convert the three arrays into sparse format, ready for porting into a linear solver. 

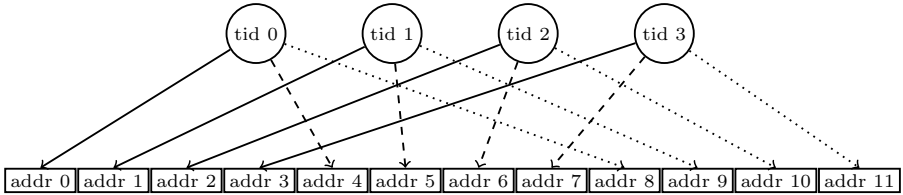
\begin{figure}[h!]
	\centering
	\begin{tikzpicture}[sharp corners=2pt,inner sep=2pt,node distance=.5cm,every text node part/.style={align=center}]
	
	\node[line width=0.25mm, circle, draw, minimum height = 2mm, minimum width = 2mm] (tid0){\scriptsize tid 0};
	\node[line width=0.25mm, circle, draw, right=1cm of tid0, minimum height = 2mm, minimum width = 2mm] (tid1){\scriptsize tid 1};
	\node[line width=0.25mm, circle, draw, right=1cm of tid1, minimum height = 2mm, minimum width = 2mm] (tid2){\scriptsize tid 2};
	\node[line width=0.25mm, circle, draw, right=1cm of tid2, minimum height = 2mm, minimum width = 2mm] (tid3){\scriptsize tid 3};
	
	\node[line width=0.25mm, draw, below left=15mm and 21mm of tid0, minimum height = 2mm, minimum width = 2mm] (addr0){\scriptsize addr 0};
	\node[line width=0.25mm, draw, right=0.01cm of addr0, minimum height = 2mm, minimum width = 2mm] (addr1){\scriptsize addr 1};
	\node[line width=0.25mm, draw, right=0.01cm of addr1, minimum height = 2mm, minimum width = 2mm] (addr2){\scriptsize addr 2};
	\node[line width=0.25mm, draw, right=0.01cm of addr2, minimum height = 2mm, minimum width = 2mm] (addr3){\scriptsize addr 3};
	\node[line width=0.25mm, draw, right=0.01cm of addr3, minimum height = 2mm, minimum width = 2mm] (addr4){\scriptsize addr 4};
	\node[line width=0.25mm, draw, right=0.01cm of addr4, minimum height = 2mm, minimum width = 2mm] (addr5){\scriptsize addr 5};
	\node[line width=0.25mm, draw, right=0.01cm of addr5, minimum height = 2mm, minimum width = 2mm] (addr6){\scriptsize addr 6};
	\node[line width=0.25mm, draw, right=0.01cm of addr6, minimum height = 2mm, minimum width = 2mm] (addr7){\scriptsize addr 7};
	\node[line width=0.25mm, draw, right=0.01cm of addr7, minimum height = 2mm, minimum width = 2mm] (addr8){\scriptsize addr 8};
	\node[line width=0.25mm, draw, right=0.01cm of addr8, minimum height = 2mm, minimum width = 2mm] (addr9){\scriptsize addr 9};
	\node[line width=0.25mm, draw, right=0.01cm of addr9, minimum height = 2mm, minimum width = 2mm] (addr10){\scriptsize addr 10};
	\node[line width=0.25mm, draw, right=0.01cm of addr10, minimum height = 2mm, minimum width = 2mm] (addr11){\scriptsize addr 11};
	
	\draw[->,line width=0.25mm] (tid0) -- (addr0.north);
	\draw[dashed,->,line width=0.25mm] (tid0) -- (addr4.north);
	\draw[dotted,->,line width=0.25mm] (tid0) -- (addr8.north);
	
	\draw[->,line width=0.25mm] (tid1) -- (addr1.north);
	\draw[dashed,->,line width=0.25mm] (tid1) -- (addr5.north);
	\draw[dotted,->,line width=0.25mm] (tid1) -- (addr9.north);
	
	\draw[->,line width=0.25mm] (tid2) -- (addr2.north);
	\draw[dashed,->,line width=0.25mm] (tid2) -- (addr6.north);
	\draw[dotted,->,line width=0.25mm] (tid2) -- (addr10.north);
	
	\draw[->,line width=0.25mm] (tid3) -- (addr3.north);
	\draw[dashed,->,line width=0.25mm] (tid3) -- (addr7.north);
	\draw[dotted,->,line width=0.25mm] (tid3) -- (addr11.north);
	\end{tikzpicture}
	\caption{Memory storage pattern of Approach 1.}
	\label{fig:coothreads}
\end{figure}

In particular, we create three arrays to store row index, column index and quadrature value. Each thread uses its own memory space on these arrays to store the computed integrals. Careful assignment of the correct memory space on each thread can drastically improve the memory store performance, taking advantage of the coalesced memory accesses. Thus, it is possible to achieve 100\% memory store efficiency. This is illustrated in Figure \ref{fig:coothreads}, where we depict the memory access pattern of each thread. Arrows with the same type (solid, dashed, or dotted lines, respectively) belong to the same iteration. Therefore, if they belong to the same warp, they write simultaneously. As we can see, threads with contiguous indices write in a coalesced fashion, which typically accelerates memory operations.

A key advantage of this first approach is that write operations can be arranged in a coalesced manner, thereby increasing performance. Moreover, since each simplex of the simplicial subdivision is regarded as a stand-alone element, mesh partitioning when using multiple GPUs (or processors, in general,) becomes immediate. At the same time, the presence of duplicates results into \Rev{fewer} elements and faces being calculated per kernel invocation. Most importantly, however, creating the final sparse matrix structure out of the three arrays becomes increasingly expensive, with the increase in the number of duplicates.

We investigate the practical performance of Approach 1, by computing the complete assembly of one time-step of the space-time stiffness matrix arising when using the dG method \eqref{spacetime_dg} on a linear parabolic problem. For each polynomial degree $p$ the mesh was chosen such that the resulting matrix could fit in the global memory of a single GPU card. On average 3 triangles are agglomerated to construct a polytopic element. All numerical investigations use a single Tesla P100 PCIe card with 16GB \Rev{of global memory}, having a total of $3584$ cores, with \Rev{a total peak double precision performance of} $4.67$ teraFLOPS. The host machine has a pair of 14-core Intel Xeon E5-2680v4 with 256GB of DDR4 memory.

\begin{table}[!h]
	\centering
	\begin{tabular}{l|
			S[table-format=1.3]
			S[table-format=2.3]
			S[table-format=2.3]
			S[table-format=2.2]
			S[table-format=2.1]}
		\hline
		& {$p=1$} & {$p=2$}  & {$p=3$}  & {$p=4$}  & {$p=5$}   \\ \hline \hline
		\multicolumn{6}{c}{\bf single precision}   \\ \hline \hline
		element kernel          & 0.007    &  0.059   &  0.27    & 1      &  3      \\ 
		interior kernel         & 0.003    &  0.024   &  0.087   & 0.29   &  0.7   \\ 
		inflow kernel           & 0.002    &  0.014   &  0.054   & 0.16   &  0.4    \\ 
		\textbf{total kernels}  & 0.014    &  0.1     &  0.41    & 1.4    &  4.1    \\ \hline
		indices                 & 4.3      &  8.6     & 18.5    & 31.8    & 49     \\ \hline \hline
		\textbf{total assembly} & 4.7      &  9.6     & 20.6    & 36      & 57     \\ \hline
		\multicolumn{6}{c}{\bf double precision}   \\ \hline \hline
		element kernel          & 0.008    &  0.067   &  0.31   & 1.2    &  3.7     \\ 
		interior kernel         & 0.005    &  0.039   &  0.14   & 0.45   &  1.1     \\ 
		inflow kernel           & 0.003    &  0.016   &  0.06   & 0.18   &  0.5     \\ 
		\textbf{total kernels}  & 0.017    &  0.12    &  0.51   & 1.8    &  5.3     \\ \hline
		indices                 & 4.7      & 10       & 18.5   & 30      & 55.2    \\ \hline \hline
		\textbf{total assembly} & 5.2      & 11.3     & 21.3   & 35      & 66.4    \\ \hline
	\end{tabular}
	\caption{Approach 1: \emph{seconds per million degrees of freedom} using single and double precision, respectively.}
	\label{tab:coo_times}
\end{table}
In Table \ref{tab:coo_times} we collect the assembly times recorded using single and double precision arithmetic, respectively. In the first three rows of each case the \emph{execution time per million degrees of freedom \Rev{(DoFs)}} for the elements, the interior faces, and the inflow faces, respectively, is presented. The total execution time for all the kernels is given, which also includes the (negligible) execution time for the boundary faces and the imposition of Dirichlet conditions. \Rev{For the CUDA kernel timings in this work, we used the \texttt{nvprof} tool and we reported the average of five invocations for each kernel.} The rows ``indices'' record the time needed to build the sparse matrix out of the three arrays (row, column, value), which includes sorting of the indices as well as addition and removal of the duplicate quadrature values. For the index sorting step, we employ the \texttt{csr\_matrix} class from Python's SciPy \cite{2020SciPy-NMeth} sparse module upon transferring all data back to the host memory from the GPU memory. Note that this built-in sorting method is single-threaded \Rev{and dominates} the total time. We prefer to include the cost of the index sorting step in the discussion, even though it is not a primal concern of the proposed work; its performance will highlight a number of important observations. Of course, it may be possible to reduce the runtime by the use of multi-threaded sorting, or investigate other options. Nevertheless, we prefer not to pursue this direction further at this point, since our primary interest lies in the assembly of polytopic meshes consisting of ``many'' simplices per element. As such, the bottle-neck of manipulation of duplicates poses a related severe challenge, which will be addressed below. 

\begin{figure}[!h]
	\centering
	\includegraphics[width=6.6cm,height=5cm]{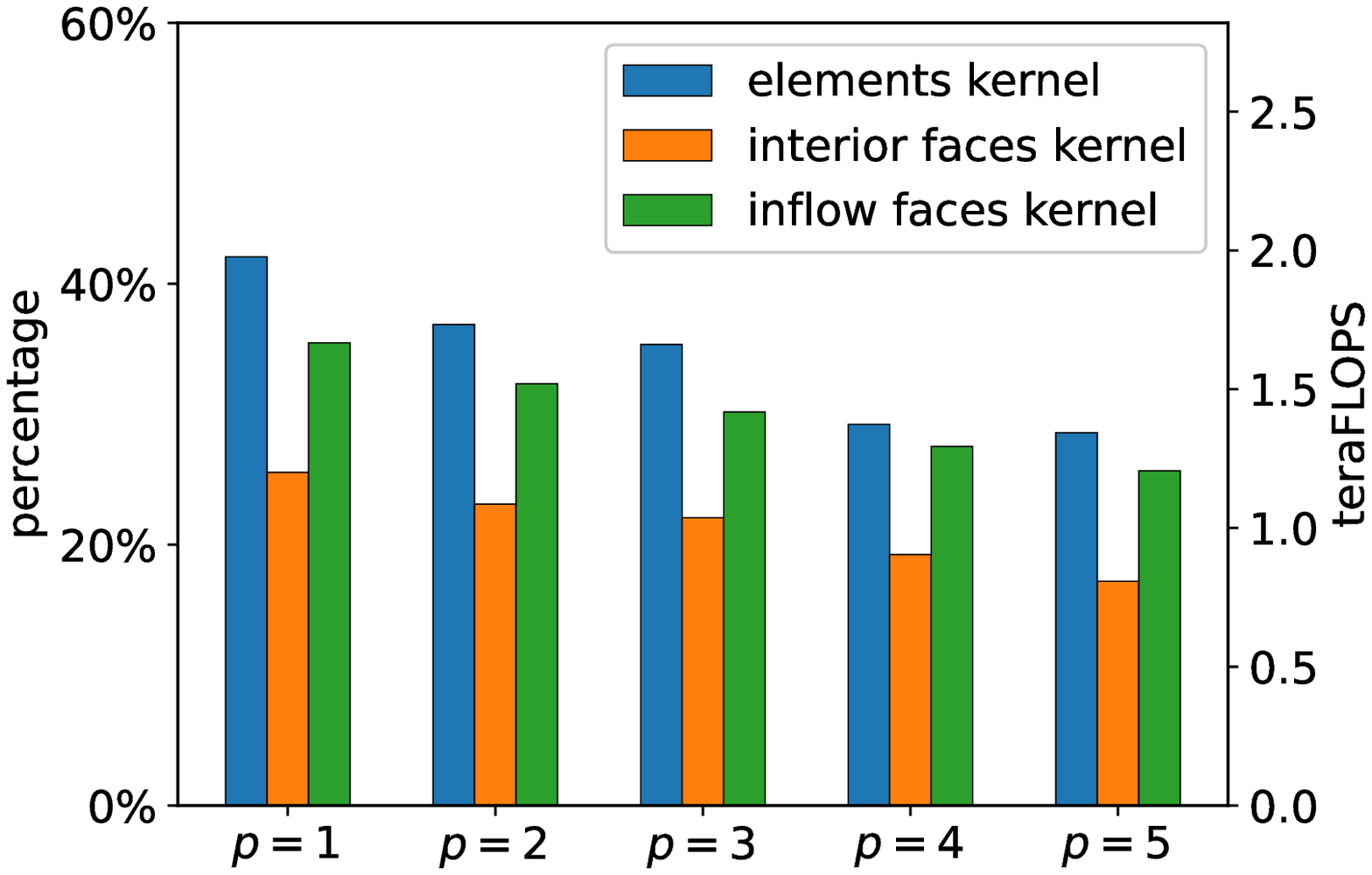}
		\includegraphics[width=6cm,height=5cm]{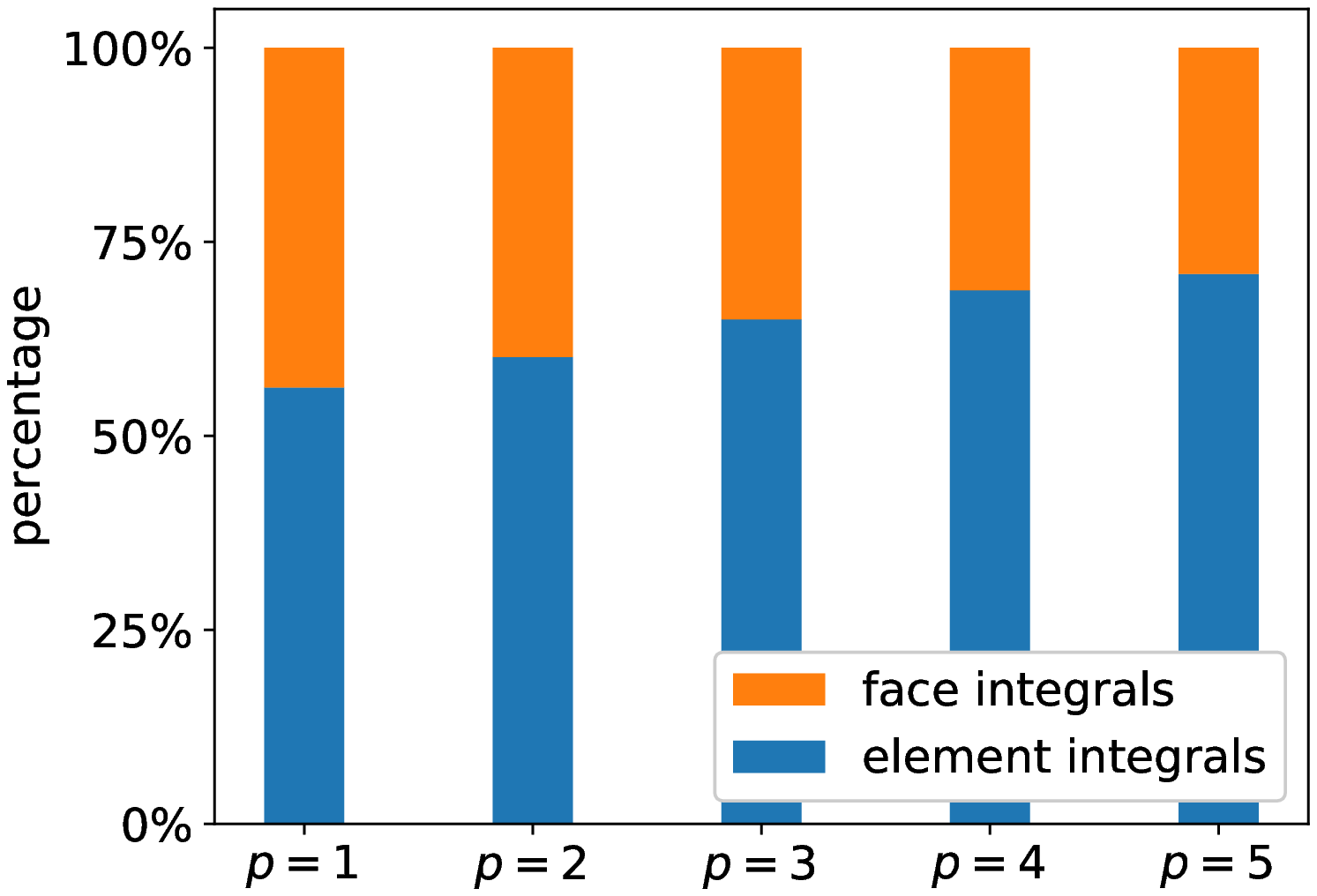}
		\caption{Approach 1: (left) double precision performance percentage and FLOPS; (right)  percentage of computation time for element vs.~face kernels, for $p=1,\dots,5$.}
	\label{fig:metrics_dp_coo}
\end{figure}

\begin{figure}[!h]
	\centering
	\includegraphics[width=6.6cm,height=5cm]{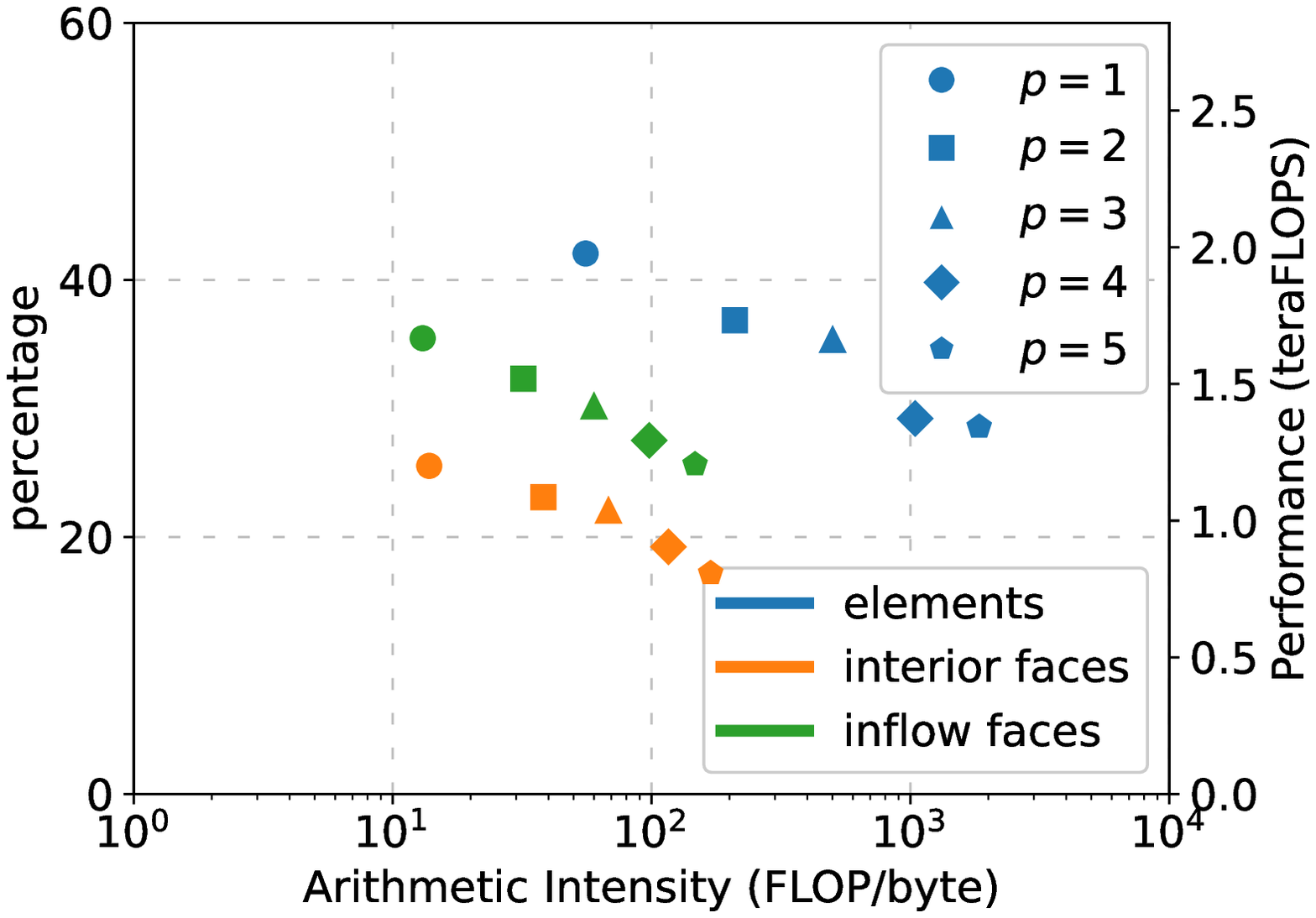}
		\includegraphics[width=6cm,height=5cm]{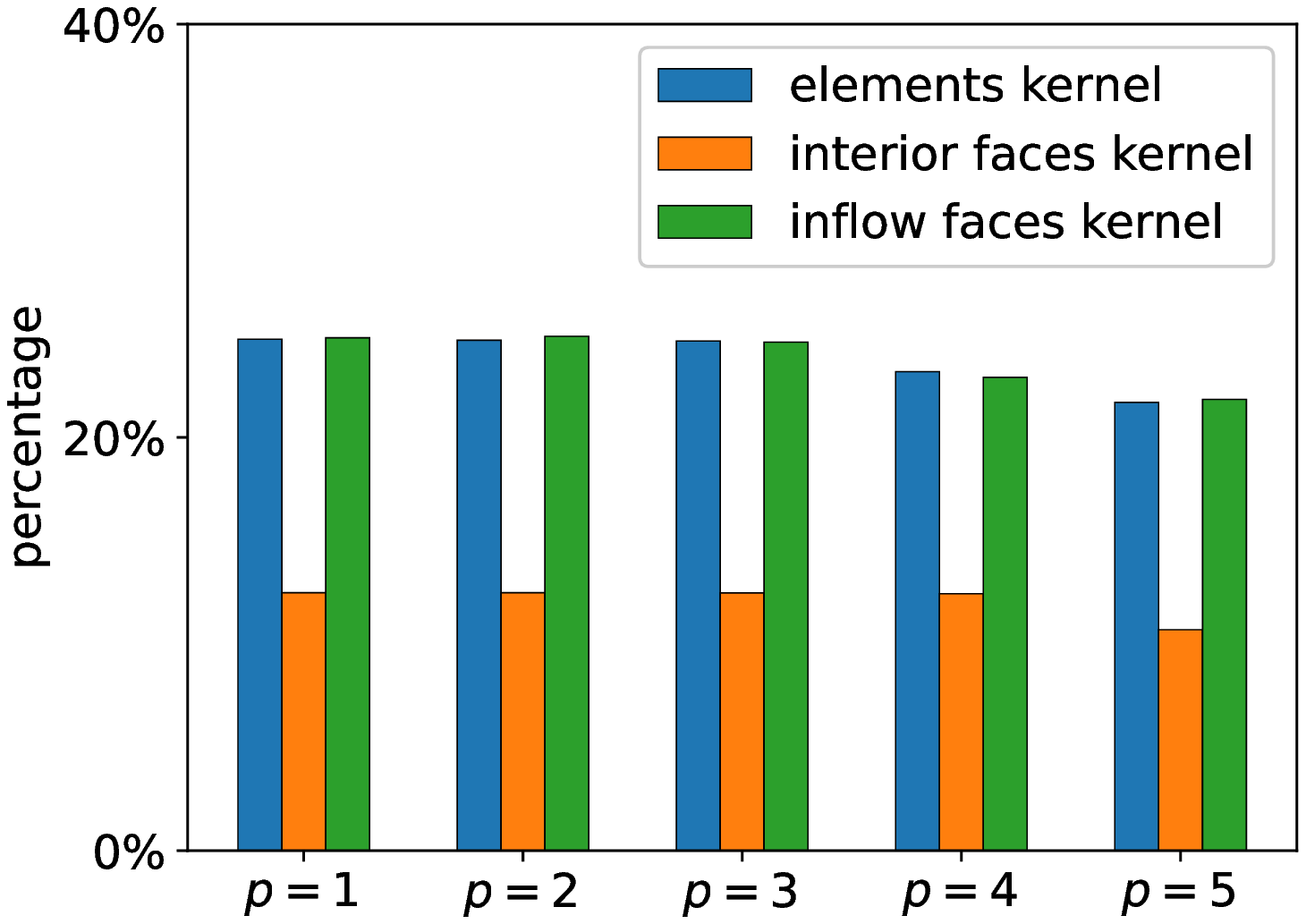}
		\caption{Approach 1: (left) performance in teraFLOPS over arithmetic intensity for transactions to/from DRAM; (right) achieved occupancy for the three main kernels.}
	\label{fig:metrics_dp_coo2}
\end{figure}

In any case, the implementation is able to assemble a million degrees of freedom using linear basis functions in less than $5$ seconds, with the integral evaluation taking only $14$ milliseconds. As the polynomial degree $p$ grows, the assembly times grow, due to the higher cost per degree of freedom by the decreasing sparsity and the increasing quadrature cost of the resulting system matrix. The ``total assembly'' is the time required to provide the complete system matrix in compressed sparse row (CSR) format, ready for the linear solver.


In Figure \ref{fig:metrics_dp_coo} (left panel), we record the performance achieved by each kernel. The maximum performance of more than $40\%$, translating into $1.9$ teraFLOPS on the Tesla P100 is achieved by the element integral kernel for $p=1$. In Figure \ref{fig:metrics_dp_coo} (right panel) we present a comparison of the total execution time for element versus the combined face integrals. The total execution time is increasingly dominated by the element integral kernel as the order of approximation increases: \Rev{each $dD$-simplicial sub-element quadrature requires $O(p^d)$ operations, while simplicial face integrals need $O(p^{d-1})$; cf. Figures \ref{fig:metrics_dp_coo} and \ref{fig:metrics_dp_csr} below.}

We note, however that, due to the \Rev{number} of duplicate values, this approach may not be recommended for meshes \Rev{comprised} of polytopic elements with \emph{many} faces per element, e.g., ones arising from agglomeration procedures; see Figure \ref{fig:30elements} below for an example of such mesh.

\subsection{Approach 2}
 In contrast, following the second approach, we compute first the non--zero indices and delete the duplicates, before storing them into a sparse format for fast access. Subsequently, we calculate the quadrature values for each simplex of the simplicial subdivision and we atomically store them into their target positions.
The key idea of the second approach is to first precompute the sparsity pattern of the stiffness and mass matrices and, subsequently, write the calculated quadrature values directly into their final position. By doing so, we can calculate more elements in the same kernel invocation, since no duplicate values are stored in the GPU's global memory. Moreover, in this way we make use of the special structure to accelerate the creation of the sparsity pattern. This is particularly relevant for polytopic meshes with \emph{many} faces per element as we shall see below.


The matrix of the resulting linear system from the dG discretization has a natural sparse block-structure. 
Given that the number of basis functions depends \emph{only} on the polynomial degree (and \emph{not} on the particular element shape) the index set for each element and face is precomputed through the knowledge of the local polynomial degree. Thus, we can avoid creating duplicate (global) indices for the elements of the simplicial subdivision. Crucially, the same principle applies to the face integral computations, whereby \emph{all} the face integrals of the common interface (containing many faces per element) between two elements are stored on the same blocks. For example, in a mesh of approximately $500k$ triangles agglomerated into $8k$ elements with $p=3$, following Approach 1 of creating all the indices for each individual interior face and converting them into CSR format took $6.4s$. In contrast, on the same processor it took $1.7s$, using Approach 2, whereby we allocate a unique index per common element interface (containing many faces).

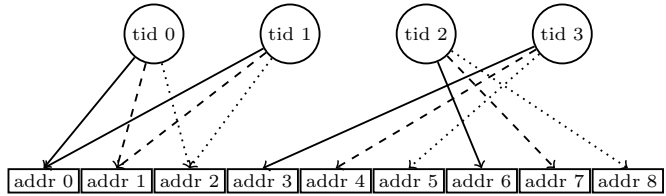
\begin{figure}[h!]
	\centering
	\begin{tikzpicture}[sharp corners=2pt,inner sep=2pt,node distance=.5cm,every text node part/.style={align=center}]
	
	\node[line width=0.25mm, circle, draw, minimum height = 2mm, minimum width = 2mm] (tid0){\scriptsize tid 0};
	\node[line width=0.25mm, circle, draw, right=1cm of tid0, minimum height = 2mm, minimum width = 2mm] (tid1){\scriptsize tid 1};
	\node[line width=0.25mm, circle, draw, right=1cm of tid1, minimum height = 2mm, minimum width = 2mm] (tid2){\scriptsize tid 2};
	\node[line width=0.25mm, circle, draw, right=1cm of tid2, minimum height = 2mm, minimum width = 2mm] (tid3){\scriptsize tid 3};
	
	\node[line width=0.25mm, draw, below left=15mm and 7mm of tid0, minimum height = 2mm, minimum width = 2mm] (addr0){\scriptsize addr 0};
	\node[line width=0.25mm, draw, right=0.01cm of addr0, minimum height = 2mm, minimum width = 2mm] (addr1){\scriptsize addr 1};
	\node[line width=0.25mm, draw, right=0.01cm of addr1, minimum height = 2mm, minimum width = 2mm] (addr2){\scriptsize addr 2};
	\node[line width=0.25mm, draw, right=0.01cm of addr2, minimum height = 2mm, minimum width = 2mm] (addr3){\scriptsize addr 3};
	\node[line width=0.25mm, draw, right=0.01cm of addr3, minimum height = 2mm, minimum width = 2mm] (addr4){\scriptsize addr 4};
	\node[line width=0.25mm, draw, right=0.01cm of addr4, minimum height = 2mm, minimum width = 2mm] (addr5){\scriptsize addr 5};
	\node[line width=0.25mm, draw, right=0.01cm of addr5, minimum height = 2mm, minimum width = 2mm] (addr6){\scriptsize addr 6};
	\node[line width=0.25mm, draw, right=0.01cm of addr6, minimum height = 2mm, minimum width = 2mm] (addr7){\scriptsize addr 7};
	\node[line width=0.25mm, draw, right=0.01cm of addr7, minimum height = 2mm, minimum width = 2mm] (addr8){\scriptsize addr 8};
	
	\draw[->,line width=0.25mm] (tid0) -- (addr0.north);
	\draw[dashed,->,line width=0.25mm] (tid0) -- (addr1.north);
	\draw[dotted,->,line width=0.25mm] (tid0) -- (addr2.north);
	
	\draw[->,line width=0.25mm] (tid1) -- (addr0.north);
	\draw[dashed,->,line width=0.25mm] (tid1) -- (addr1.north);
	\draw[dotted,->,line width=0.25mm] (tid1) -- (addr2.north);
	
	\draw[->,line width=0.25mm] (tid2) -- (addr6.north);
	\draw[dashed,->,line width=0.25mm] (tid2) -- (addr7.north);
	\draw[dotted,->,line width=0.25mm] (tid2) -- (addr8.north);
	
	\draw[->,line width=0.25mm] (tid3) -- (addr3.north);
	\draw[dashed,->,line width=0.25mm] (tid3) -- (addr4.north);
	\draw[dotted,->,line width=0.25mm] (tid3) -- (addr5.north);
	\end{tikzpicture}
	\caption{Memory storage pattern of Approach 2 with scattered memory operations. Lines with same type are write operations that are performed simultaneously. Threads 0 and 1 are calculating integrals of the same polygonal element, hence they write in the same memory locations with atomic operations.} 
	\label{fig:csrthreads}
\end{figure}


\begin{table}[!h]
	\centering
	\begin{tabular}{l|
			S[table-format=1.3]
			S[table-format=1.3]
			S[table-format=1.2]
			S[table-format=1.2]
			S[table-format=2.1]}
		\hline
		& {$p=1$} & {$p=2$}  & {$p=3$}  & {$p=4$}  & {$p=5$}   \\ \hline \hline
		\multicolumn{6}{c}{\bf single precision}   \\ \hline \hline
		element kernel          & 0.007   & 0.061   & 0.27  & 1     &  3     \\ 
		interior kernel         & 0.019   & 0.09    & 0.2   & 0.36  &  0.8   \\ 
		inflow kernel           & 0.002   & 0.015   & 0.06  & 0.16  &  0.4   \\ 
		\textbf{total kernels}  & 0.029   & 0.17    & 0.52  & 1.5   &  4.2   \\ \hline
		indices                 & 0.63    & 1.2     & 2.3   & 3.9   &  7     \\ \hline \hline
		\textbf{total assembly} & 0.79    & 1.5     & 3.1   & 6     & 12.2   \\ \hline
		\multicolumn{6}{c}{\bf double precision}   \\ \hline \hline
		element kernel          & 0.008   & 0.069   & 0.3   & 1.2   &  3.7   \\ 
		interior kernel         & 0.018   & 0.083   & 0.2   & 0.53  &  1.2   \\ 
		inflow kernel           & 0.003   & 0.017   & 0.07  & 0.19  &  0.5   \\ 
		\textbf{total kernels}  & 0.03    & 0.17    & 0.57  & 1.9   &  5.5   \\ \hline
		indices                 & 0.84    & 1.5     & 2.7   & 4.4   &  7.7   \\ \hline \hline
		\textbf{total assembly} & 1.02    & 1.9     & 3.7   & 7.2   & 14.7   \\ \hline
	\end{tabular}
	\caption{Approach 2: seconds per million degrees of freedom using single and double precision, respectively.}
	\label{tab:csr_times}
\end{table}


\begin{figure}[!h]
	\centering
	\includegraphics[width=6.6cm,height=5cm]{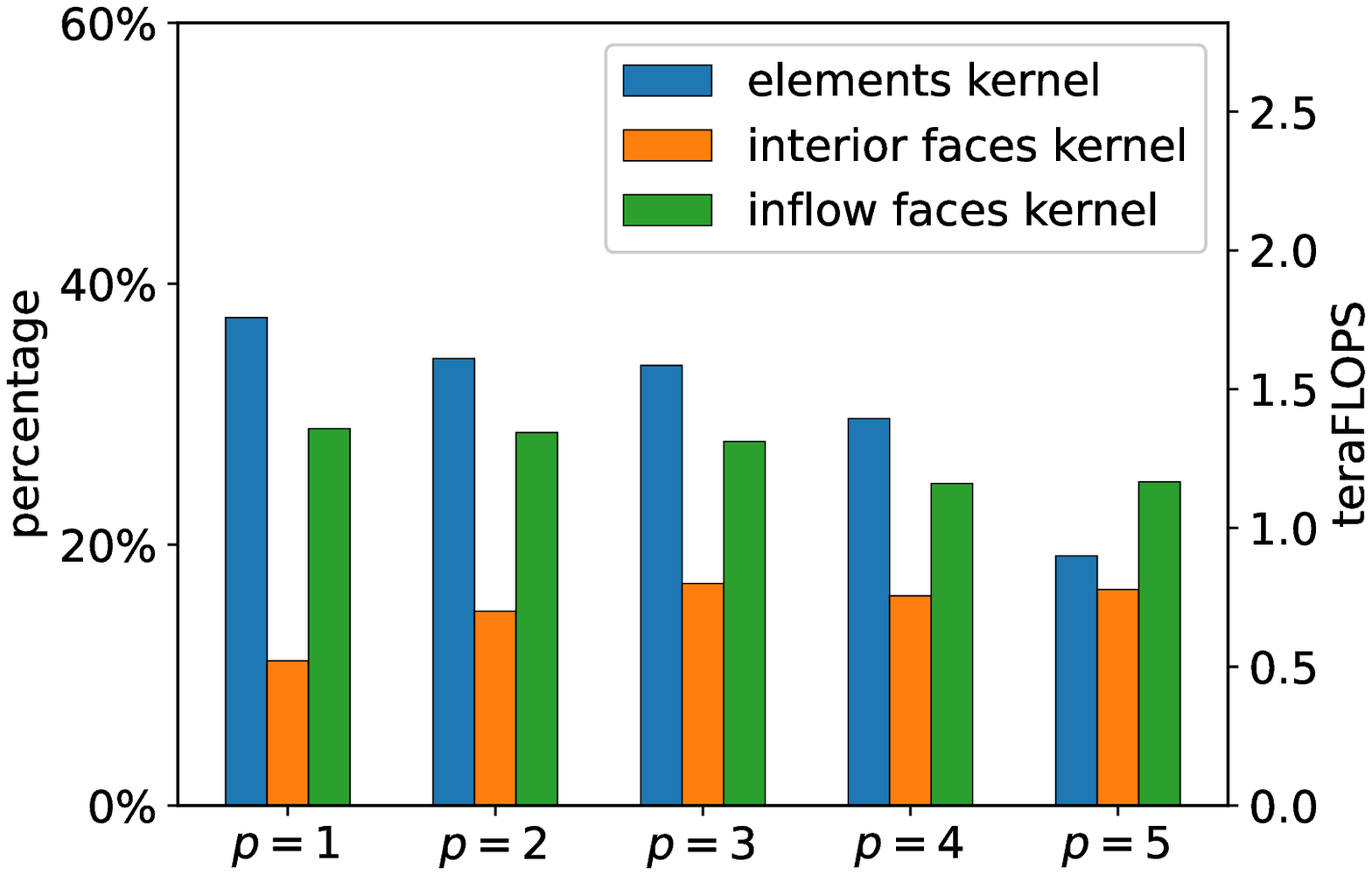}	\includegraphics[width=6cm,height=5cm]{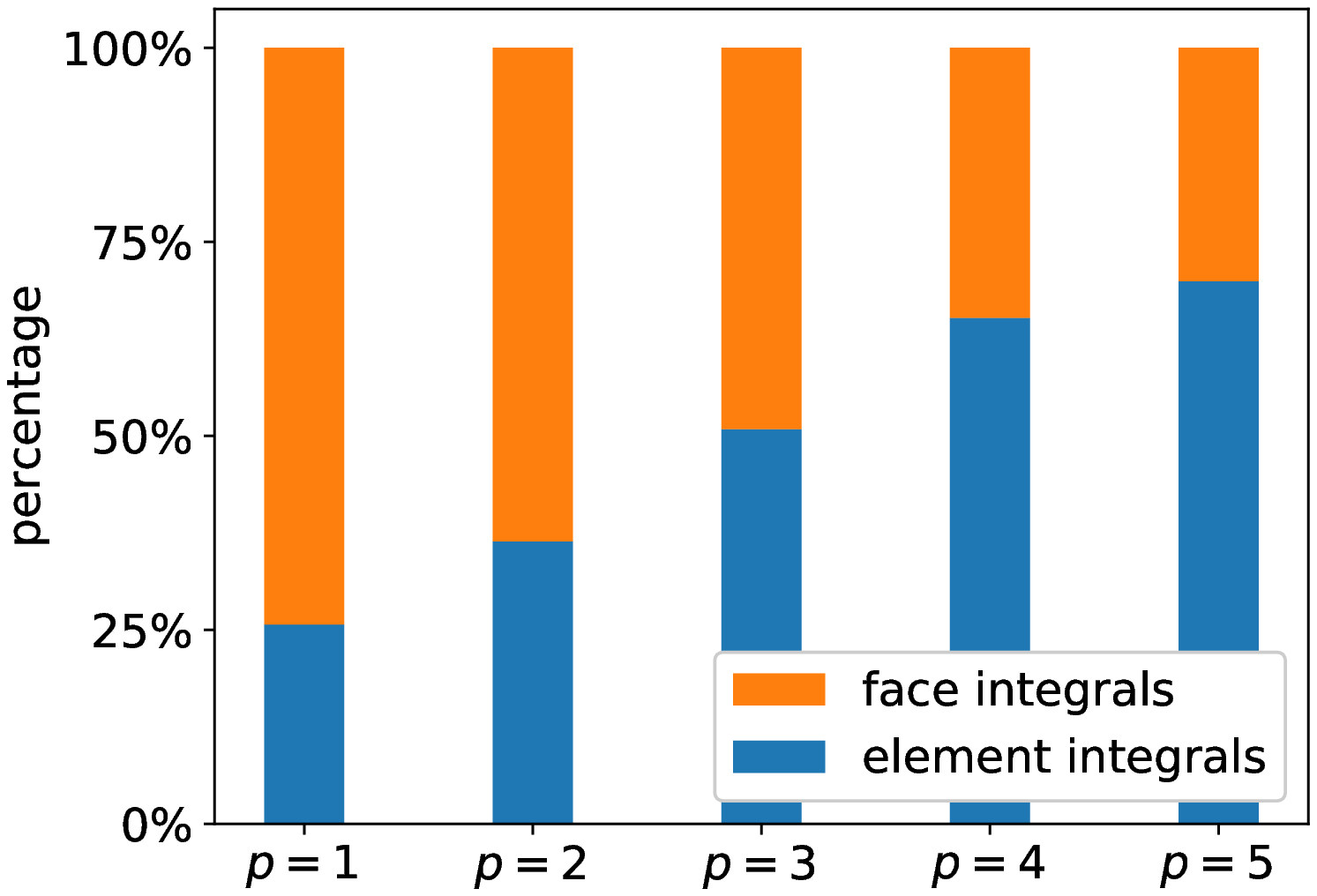}
	\caption{Approach 2: (left) double precision performance percentage and FLOPS; (right)  percentage of computation time for element vs.~face kernels, for $p=1,\dots,5$.}
	\label{fig:metrics_dp_csr}
\end{figure}

\begin{figure}[!h]
	\centering
	\includegraphics[width=6.6cm,height=5cm]{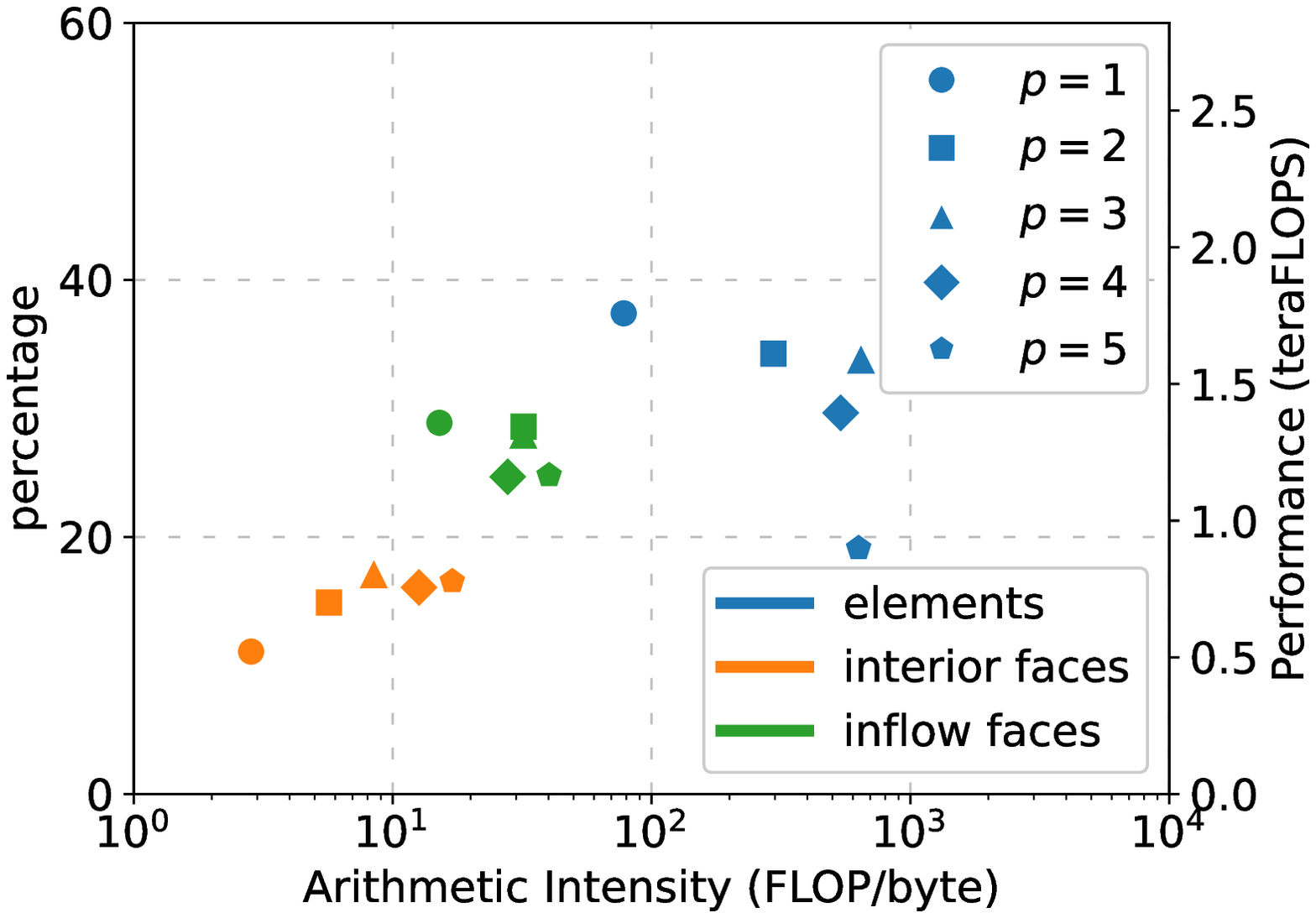}	\includegraphics[width=6cm,height=5cm]{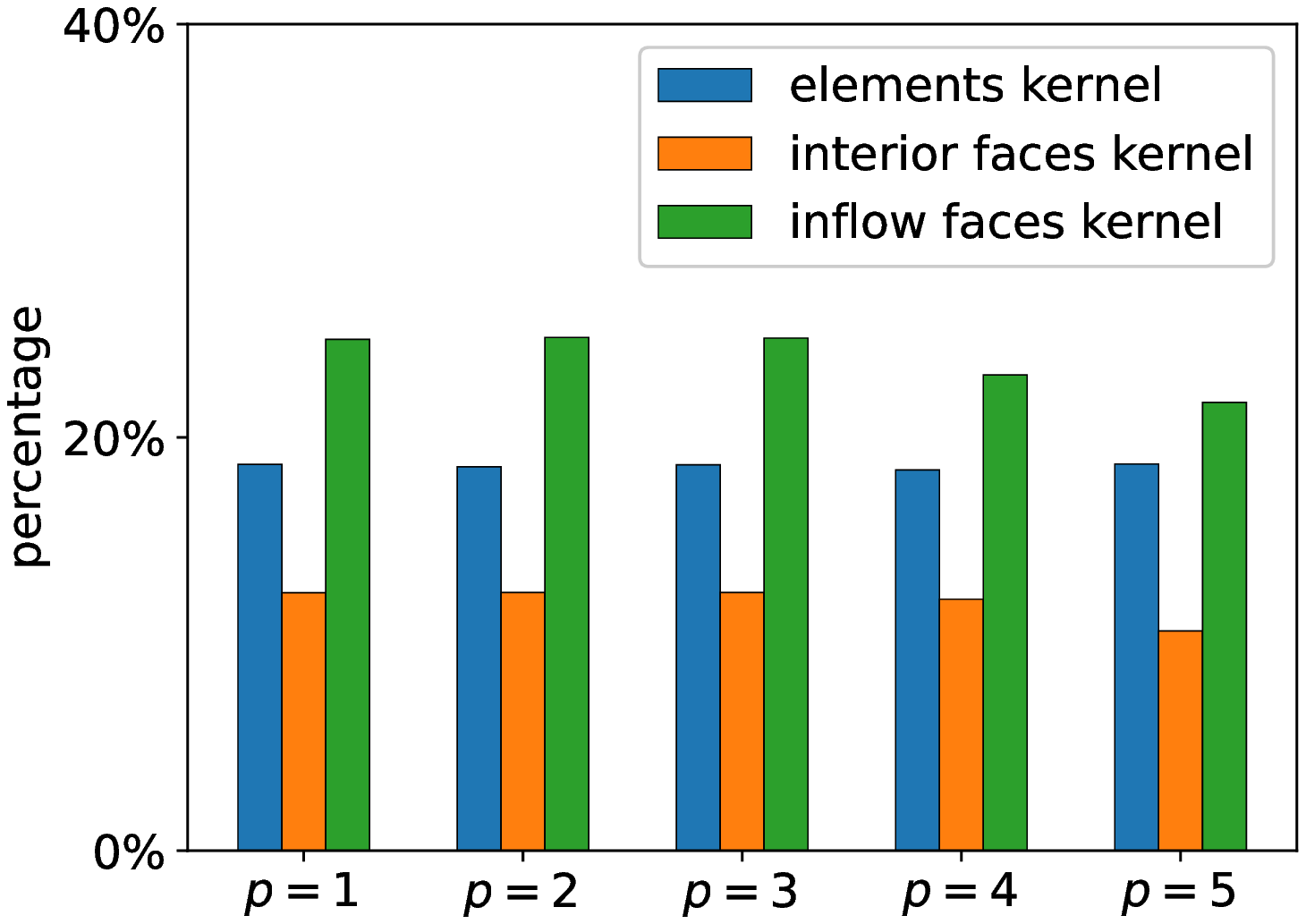}
	\caption{Approach 2: (left) performance in teraFLOPS over arithmetic intensity for transactions to/from DRAM; (right) achieved occupancy for the three main kernels.}
	\label{fig:metrics_dp_csr2}
\end{figure}

Upon creation, the matrix is transferred to the device memory for population with the calculated quadrature values. To identify the position where each value must be added, a binary search is required. Care must be taken when adding the contributions of every simplex or face to rule out, so-called, race conditions. Since it is possible that simplices belonging to the same polygon can also belong to the same execution warp, they will try to update the same memory locations in global memory simultaneously, as is illustrated in Figure \ref{fig:csrthreads}. To avoid this, these writes must be made as atomic operations. It is worth mentioning here that since floating point addition is not associative and there is no a priori guarantee on the order of the atomic additions, different matrices (up to floating point precision) may be created after multiple executions of the same exact problem. Nonetheless, in the numerous numerical investigations we performed, this issue \Rev{appears to have negligible effect on} the results.


In Table \ref{tab:csr_times} we collect the assembly times recorded using single and double precision arithmetic, respectively. Here, \Rev{the ``indices'' row records the time spent on creating, sorting and converting the indices into the CSR format}, which now precedes the execution of the kernels. With this approach, the time to assemble a million degrees of freedom, using linear basis functions is less than $1$ second, almost six times faster than Approach 1. In Figure \ref{fig:metrics_dp_csr} (left panel), we can see the performance achieved by each kernel. The scattered memory access pattern affects negatively the performance, mainly for low order methods, compared to Approach 1. For high order methods calculating the quadrature value becomes more expensive compared to the memory operations to store it in its final position, so the scattered memory access pattern has a marginal effect. Finally, in Figure \ref{fig:metrics_dp_csr} (right panel), a comparison of the total execution time for element versus face integrals is recorded. The face integral kernel dominates the \Rev{total execution time for low order elements,} resulting to a modest increase in the total kernel execution time, compared to Approach 1. Nonetheless, the significant savings recorded in the index manipulation of Approach 2, (cf.~Tables \ref{tab:coo_times} and \ref{tab:csr_times},) showcases its viability. We note that both approaches use the same index sorting functions of Python's SciPy module. More sophisticated/parallel implementation of index sorting routines has the potential to improve the ``indices'' (and, of course, ``total'') times for both approaches. Indeed, as we shall see in Section \ref{MPI} below, the index sorting run time can be reduced via an implementation on multiple GPUs.

\begin{table}[!h]
	\centering
	\begin{tabular}{ l | c c }
		\hline
		 & Approach 1 & Approach 2 \\
		\hline \hline
		elements & 126 & 144  \\
		interior faces & 198 & 210 \\
		inflow faces & 109 & 107 \\
		\hline
	\end{tabular}
	\caption{Number of 32bit registers used per thread, as reported by \texttt{nvprof}.}
	\label{tab:registers_per_thread}
\end{table}

\begin{remark}
\Rev{The arithmetic intensity results in Figures \ref{fig:metrics_dp_coo2} and \ref{fig:metrics_dp_csr2} (left panels) suggest that we are, in most cases, in the compute bound region. The theoretical occupancy for all kernels, based on the selected threads per block and used registers per thread, matches the achieved occupancy, cf. right panels in Figures \ref{fig:metrics_dp_coo2} and \ref{fig:metrics_dp_csr2}. While it would be possible to reduce the number of registers per thread in an effort to further increase performance, such a direction would limit considerably the sought-after generality of the kernels, which is central to this work: the kernels have been designed to handle essentially arbitrary polynomial degrees and element shapes.}

\end{remark}

\section{Numerical experiments}\label{sec:numerics}
We continue the investigation of the performance of the algorithms by testing them on two \Rev{computational  problems.} The first deals with the question of performance on meshes comprising elements with \emph{many} faces each, constructed by  aggressive agglomeration of a fine background mesh. The second series of numerical experiments investigates the performance and scalability of an implementation of the algorithms on multiple GPUs. 

\subsection{Performance on highly agglomerated meshes}
Highly agglomerated meshes, i.e., meshes  arising from agglomerating \emph{many} simplices are relevant in many applications \cite{DGpolybook}, such as domains with highly heterogeneous boundaries \cite{DGease,MR2846986} and also in the context of multilevel solvers \cite{MR3585793,hp_multigrid_polytopes_2017}. Approach 1 is not suitable in this setting due to the excessive number of duplicate entries. Therefore, we seek to test whether Approach 2 has a performance penalty when the mesh consists of highly agglomerated elements, due to the high number of atomic operation replays arising by large numbers of threads updating the same memory locations simultaneously. Nevertheless, the benefit from not creating duplicate indices of Approach 2 in this context is expected to offer superior overall performance of the assembly process. 

\begin{figure}[!h]
	\centering
	\includegraphics[scale=0.5]{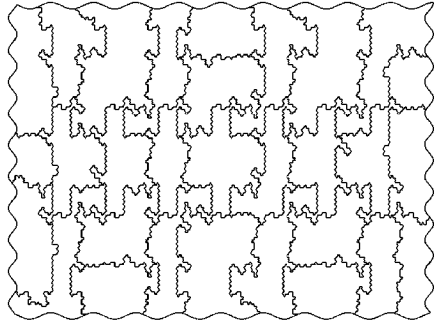}
	\caption{$503,596$ triangles agglomerated into $30$ polygonal elements.}
	\label{fig:30elements}
\end{figure}

\Rev{We start with a problem defined on the domain \Rev{$\Omega$ with oscillating boundary, which is approximating $(0,1)^2$. To} represent the computational domain, we employed $503,596$ unstructured simplicial elements as the background mesh}; we refer to Figure \ref{fig:30elements} for an illustration of the extreme agglomeration process resulting to a $30$-element polygonal mesh. \Rev{We note that exactly the same dG method on exactly these meshes have been used in \cite{DGease} for the numerical approximation of a convection-diffusion problem, where optimal convergence rates have been observed and recorded.} In our computations below, the original simplicial mesh is agglomerated into $8,337$ and $125,981$ elements respectively. We also test the same method by assembling on the original mesh of $503,596$ simplices, treating each simplex as an element. The characteristics of the meshes can be found in Table \ref{tab:agglomerated_details}. The resulting $2D$ polygonal meshes are tensorised into forming $3D$ prismatic space-time elements. On the latter, we assemble the space-time dG method \eqref{spacetime_dg} with the reduced-complexity $\mathcal{P}_{p_{\kappa_n}}$ Galerkin space choice, for the numerical approximation of the linear parabolic problem \eqref{parabolic-problem1} with the specifics: 
$
{\bf a} = I_2$, ${\bf w}= (0,0)^\top$, $c=1$.
The load function $f$ is selected so that the solution $u(x,y,t)=\sin{(\pi x)} \sin{(\pi y)} (1-t)$.

\begin{table}[!h]
	\centering
	\begin{tabular}{ l | c c c }
		\hline
		& Mesh 1 & Mesh 2 & Mesh 3 \\
		\hline \hline
		\#elements & 8,337 & 125,981 & 503,596 \\
		\#triangles & 503,596 & 503,596 & 503,596 \\
		\#interior faces & 124,628 & 376,455 & 754,118 \\
		\hline
	\end{tabular}
	\caption{The original triangular mesh with $ 503,596$ triangles, and meshes constructed via two levels of agglomeration of the original mesh.}
	\label{tab:agglomerated_details}
\end{table}

\begin{figure}[!h]
	\centering
	\includegraphics[scale=0.55]{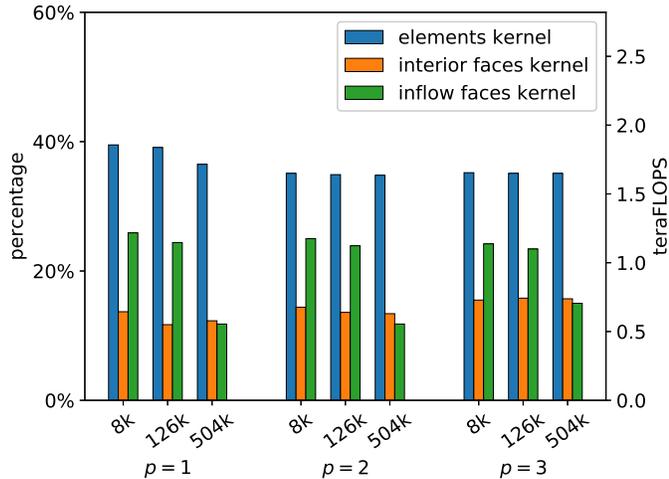}
	\caption{Double precision performance of the three main kernels on the $3$ meshes from Table \ref{tab:agglomerated_details}, for $p=1,2,3$.}
	\label{fig:dp_performance_agglomerated}
\end{figure}

For this test we used a single Tesla P100 PCIe card with 16GB of 4096 bit HBM2 \Rev{global} memory, with a total of 3584 cores and processing power of 4.67 teraFLOPS.
In Figure \ref{fig:dp_performance_agglomerated}, we record the performance achieved by Approach 2 in double precision. Similar performance is observed in all cases. The element integrals kernel achieved between 35\% and 40\% of peak performance with a maximum of about $1.85$ teraFLOPS. The interior faces kernel performance was between 12\% to 16\% of the peak with a maximum computational throughput of $750$ gigaFLOPS. We also observe a performance penalty when using the original $504k$-element mesh in the inflow faces kernel. This is due to the scattered memory reads of the solution from the previous time-step. This overhead does not appear for the two meshes with $8,337$ and $125,981$ agglomerated polygonal prisms. This is due to contiguous threads that calculate quadratures from simplices belonging to the same element, reading the same solution coefficients from the previous time-step.

\begin{table}[!h]
	\centering
	\begin{tabular}{l|
			S[table-format=1.3]
			S[table-format=1.2]
			S[table-format=2.1]}
		\hline
		& {$p=1$} & {$p=2$}  & {$p=3$}   \\ \hline \hline
		\multicolumn{4}{c}{{\bf Mesh 1:} $8k$ elements}   \\ \hline \hline
		total degrees of freedom	& {$33k$} & {$83k$} & {$167k$} \\ \hline 
		element kernel          & 0.006   & 0.12     & 1      \\ 
		interior kernel         & 0.003   & 0.04     & 0.3    \\ 
		inflow kernel           & 0.002   & 0.03     & 0.2    \\ 
		\textbf{total kernels}  & 0.012   & 0.2      & 1.5    \\ \hline
		indices                 & 0.48    & 0.7      & 1.3    \\ \hline \hline
		\textbf{total assembly} & 1.34    & 1.7      & 3.7    \\ \hline
		
		\multicolumn{4}{c}{{\bf Mesh 2:}  $126k$ elements}   \\ \hline \hline
		total degrees of freedom	& {$504k$} & {$1.3m$} & {$2.5m$} \\ \hline 
		element kernel          & 0.006   & 0.12     & 1      \\ 
		interior kernel         & 0.012   & 0.14     & 0.9    \\ 
		inflow kernel           & 0.002   & 0.03     & 0.2    \\ 
		\textbf{total kernels}  & 0.02    & 0.3      & 2.1    \\ \hline
		indices                 & 1.15    & 3.2      & 10.7   \\ \hline \hline
		\textbf{total assembly} & 1.9     & 4.6      & 15     \\ \hline
		
		\multicolumn{4}{c}{{\bf Mesh 3:} $504k$ elements}   \\ \hline \hline
		total degrees of freedom	& {$2m$} & {$5m$} & {$10m$} \\ \hline 
		element kernel          & 0.006   & 0.12     & 1      \\ 
		interior kernel         & 0.022   & 0.3      & 1.8    \\ 
		inflow kernel           & 0.005   & 0.05     & 0.3    \\ 
		\textbf{total kernels}  & 0.033   & 0.5      & 3.1    \\ \hline
		indices                 & 2       & 7.4      & 25.7   \\ \hline \hline
		\textbf{total assembly} & 2.9     & 9.6      & 32.9   \\ \hline
	\end{tabular}
	\caption{Approach 2: seconds using double precision arithmetic for the highly agglomerated meshes and for one time-step.}
	\label{tab:csr_times_highly_agglomerated}
\end{table}

In Table \ref{tab:csr_times_highly_agglomerated} we record the \emph{kernel} and \emph{total} assembly times for $p=1,2, 3$ using double precision. Although it is not necessary to use the very fine background mesh for the quadrature evaluation (it is enough to use any sub-triangulation into as many simplices as faces), we do so in this numerical experiment to highlight further the acceleration potential of the proposed approach. This results to the execution time for the element kernels to be the same in all cases. As expected, differences in performance arise from the interior faces kernel, as the highly agglomerated meshes have far \Rev{fewer} interior faces than the background mesh. Of course, more aggressive agglomeration results to \Rev{coarser meshes and therefore fewer global degrees of freedom}. We envisage to apply such meshes within a mesh-adaptive Galerkin framework, thereby equilibrating the local resolution requirements for a given accuracy with the computational cost.

\subsection{Performance on multiple GPUs}\label{MPI}
To assess the performance and scalability of the proposed algorithms on larger scale problems, we consider a 3D problem with non-negative characteristic form \eqref{CDR} with coefficients $A=0.01 I_3$, ${\bf b}=(
1+x, 1+y,1+z)^\top$, $c=3+xyz$,
with $\Omega=(0,1)^3$. The load function $f$ is selected so that the solution $u(x,y,z)=\sin{(\pi x)} \sin{(\pi y)} \sin{(\pi z)}$. Note that different choices of $A$ and ${\bf b}$ result to different sets of active face terms on each element; the above choice is a typical scenario.  The polyhedral elements stem from the agglomeration of a fine three-dimensional unstructured tetrahedral mesh. This is in contrast to the prismatic meshes used for the parabolic problem above, as this is now a fully $3D$ unstructured grid. As such, we expect both more kernel evaluations (more terms in the bilinear form) and higher connectivity (more non-zero entries).

To estimate the scalability of the assembly process, we implement Approach 2 on clusters consisting of multiple GPUs per node. A Message Passing Interface (MPI) implementation distributes the load to each GPU card and each GPU is responsible for computing a part of the global matrix. Specifically, the following processes are implemented:

\begin{enumerate}
	\item METIS \cite{MR1639073} is used to subdivide the mesh in $Ngpus$ parts. 
	 We assign to each polyhedron a weight to minimize the communication cost and to simultaneously balance the quadrature cost among the GPUs;
	
	\item we flag the interior faces on the boundaries of the partition created by METIS. Those faces are processed 
	with a modified interior faces kernel
	;
	
	\item each GPU creates the sparsity pattern of its allocated subdivision. Instead of SciPy's built-in CPU routines, (as done in the single GPU examples above,) we use the CUDA Thrust library \cite{BELL2012359} to perform the index manipulation steps directly on the GPUs;
	
	\item 
	assembly of the partial stiffness matrices takes place on each GPU by executing the quadrature evaluation kernels.
\end{enumerate}

%
%

The MPI implementation of Approach 2 has two significant benefits. First, each GPU creates only a partial matrix in CSR format, allowing for much larger problems to be assembled in the same runtime. Also this allows to build in parallel stiffness matrices that are too large to fit in the global memory of one single card. Moreover, as index sorting is performed separately for each (smaller) partial matrix, the index computation cost is reduced. The partial matrices can then be used for local matrix--vector product operations on each GPU before communicating partial solutions with cards holding neighbouring subdivisions. This is particularly pertinent in the context of multilevel/domain decomposition algorithms.


\begin{table}[!h]
	\centering
	\begin{tabular}{l|
			S[table-format=2.1]
			S[table-format=2.1]
			S[table-format=2.1]
			S[table-format=2.1]
			S[table-format=1.1]}
		\hline
		{\#GPUs}	& {1} & {2}  & {4}  & {8}  & {16}   \\ \hline \hline
		\multicolumn{6}{c}{{\bf $p=1$}, $1.57m$ elements ($6.28m$ DoFs), $6.29m$ tetrahedra} \\ \hline \hline
		index manipulation      &  7.4    & 4.5    & 2.7   &  1.6   &  1.3    \\ 
		quadrature evaluation   &  2      & 1.7    & 1.7   &  1.5   &  1.4    \\ \hline
		\textbf{total assembly} &  9.4    & 6.2    & 4.4   &  3.1   &  2.7    \\ \hline
		\multicolumn{6}{c}{{\bf $p=2$}, $1.57m$ elements ($15.7m$ DoFs), $6.29m$ tetrahedra} \\ \hline \hline
		index manipulation      & 28.5    & 14.5   &  8.3  &  4.5   &  2.7    \\ 
		quadrature evaluation   &  9.7    &  5     &  3.9  &  2.5   &  2.3    \\ \hline
		\textbf{total assembly} & 38.2    & 19.5   & 12.2  &  7     &  5      \\ \hline
		\multicolumn{6}{c}{{\bf $p=3$}, $197k$ elements ($3.94m$ DoFs), $786k$ tetrahedra}   \\ \hline \hline
		index manipulation      & 15.1    &  7.7   &  4.7  &  2.7   &  1.6    \\ 
		quadrature evaluation   &  5.7    &  3.5   &  2.6  &  2.1   &  1.9    \\ \hline
		\textbf{total assembly} & 20.8    & 11.2   &  7.3  &  4.8   &  3.5    \\ \hline
	\end{tabular}
	\caption{Time in seconds for the assembly of the fully $3D$ problem using Approach 2 and double precision.}
	\label{tab:3d_times}
\end{table}

\Rev{The quadrature times recorded here include the kernels' execution time (as per table in the previous section) \emph{and} also the memory transfers from RAM to the global memory of GPUs.} This is done in an effort to showcase realistic assembly times. The results showcase reduction in both index manipulation and quadrature evaluation as the number of GPUs is increased. \Rev{Timings for the MPI code were performed using Python's \texttt{perf\_counter()} from the \texttt{time} module, upon calling an \texttt{MPI\_Barrier()}. Therefore, each reported time is the maximum across all MPI processes.} An interesting observation is that index manipulation run times using standard, freely available algorithms can be balanced with quadrature kernel execution times in multiple GPU architectures.

\section*{Acknowledgements}
We gratefully acknowledge the availability of the ALICE High Performance Computing Facility at the University of Leicester and also resources by the Cambridge Service for Data Driven Discovery (CSD3) operated by the University of Cambridge Research Computing Service (www.csd3.cam.ac.uk), provided by Dell EMC and Intel using Tier-2 EPSRC funding (capital grant EP/P020259/1), and DiRAC STFC funding (www.dirac.ac.uk). 
%
%
%
%



\appendix
\section{The discontinuity-penalization parameter}\label{penalty}  For completeness, we now give a precise formula for the discontinuity-penalization function $\sigma$ appearing in the IP-dG formulations \eqref{adv_galerkin_dg} and \eqref{spacetime_dg}.  The stability and error analysis of the IP-dG method under this choice of penalization, as well as a detailed discussion on the practical relevant of this choice, can be found in \cite{DGpolybook,DGease}.

Obviously, each polytopic element $\kappa\in \mathcal{T}_h$ (or $\kappa\in \mathcal{D}_h$, respectively,) can be covered by different families of simplices $\mathcal{K}_{\kappa}:=\{K_j\}_{j=1}^{m_\kappa^{}}$, of possibly different cardinalities, i.e., we have $\kappa\subset \cup_{K_j\in \mathcal{K}_{\kappa}} K_j$.  Each such family $\mathcal{K}_{\kappa}$ will be referred to as a \emph{covering} of $\kappa$ and we shall denote by $\mathbb{K}_{\kappa}$ the set of all such coverings. For instance, any subtriangulation of $\kappa$ is a valid covering; equally coverings with overlapping simplices are also valid. 
Let $h_\omega$, $\rho_\omega$  and $|\omega|$ denote the diameter, the inscribed radius, and the $sD$-volume of a domain $\omega\subset\mathbb{R}^s$, $s=1,\dots, d$, respectively. We say that an element $\kappa$ is \emph{$p$-coverable} if there exists at least one covering $\mathcal{K}_{\kappa}$ of $\kappa$, such
that: 1) $h_{K_j}\sim h_{\kappa}$ $\rho_{K_j}\sim \rho_{\kappa}$, \emph{and} 2) $\max_{{\bf x}\in \partial\kappa, {\bf z}\in \partial K_j}|{\bf x}-{\bf z}|\le \rho_{K_j}/(8p)^2$, for all $j=1,\dots,m_{\kappa}^{}$, with $|\cdot|$ denoting the Euclidean distance. In other words, an element is $p$-coverable if there exists a covering $\mathcal{K}_\kappa$ comprising simplices each with similar shape-regularity to the original element, which cover $\kappa$ within a distance at most $\rho_{K_j}/(8p)^2$ each, away from the element's boundary. In \cite{DGease}, a considerably weaker concept of $p$-coverability is used allowing, in particular, $K_j$ to be general \emph{curved} prisms. 

Now let $F\in \mathcal{F}_h^I$ a face shared by two elements $\kappa_1,\kappa_2\in\mathcal{T}_h$; if $F\subset \partial\Omega_{\rm D}$, we set $\kappa_2=\emptyset$. Denote also by $K^F\in\underline{\mathbb{K}}_\kappa$ a sub-simplex having $F$ as face also. We define the discontinuity-penalization parameter on $F$ by
\[
\sigma|_F:=C_{\sigma} \max_{\kappa\in\{\kappa_1,\kappa_2\}}\bigg\{ \min \Big\{\frac{|\kappa|}{\sup_{K^F\in\underline{\mathbb{K}}_\kappa}
	|K^F|},C_{\rm cov}(\kappa) \Big\}
	\frac{\bar{a}_\kappa p_\kappa^2|F|}{|\kappa|} \bigg\} ,
\]
for a computable constant $C_{\sigma}>0$, with $\bar{a}_\kappa:=\|{\bf n}^\top{\bf A}{\bf n}\|_{L_\infty(\kappa)}$, (correspondingly $\bar{a}_\kappa:=\|{\bf n}^\top{\bf a}{\bf n}\|_{L_\infty(\kappa)}$), and $C_{\rm cov}(\kappa):=p_\kappa^{2(d-1)}$ if $\kappa$ is $p$-coverable, or $C_{\rm cov}(\kappa):=\infty$ if not.

\bibliographystyle{siamplain}
\bibliography{Reference}
\end{document}